 \documentclass[a4paper,10pt]{article}
 \usepackage[authoryear,round]{natbib}
\usepackage{a4wide}
\usepackage[tbtags,sumlimits,intlimits,namelimits,reqno,fleqn]{amsmath}
\usepackage{amsthm}
\usepackage{amsfonts}
\usepackage{latexsym}
\usepackage{amstext}
\usepackage{amssymb}
\usepackage{amsxtra}
\usepackage{longtable}
\usepackage{authblk}
\usepackage{ifthen}   
\usepackage[english]{babel}
\usepackage{enumerate}
\usepackage{enumitem}   
\usepackage{comment}
\usepackage[pdfpagemode=UseOutlines,urlcolor=blue,citecolor=darkgreen,colorlinks=true]{hyperref}
\usepackage[pdftex]{graphicx}
\usepackage[pdftex]{xcolor}
\definecolor{darkgreen}{rgb}{0,0.35,0}
\hypersetup{citecolor=blue}
\usepackage{bibentry}
\usepackage{natbib}
\usepackage{multirow}
\usepackage{relsize}
\usepackage{tikz}
\usetikzlibrary{patterns}

\bibpunct{(}{)}{;}{a}{}{,}

\numberwithin{equation}{section}
\newcommand{\revised}[1]{{\color{black}#1}}
\newcommand{\corrected}[1]{{\color{black}\hypersetup{linkcolor=black,citecolor=black,urlcolor=black}#1}}

\newtheorem{Prop}{Proposition}[section]

\newtheorem{lem}[Prop]{Lemma}

\newtheorem{Th}[Prop]{Theorem}

\newtheorem{Remark}[Prop]{Remark}

\newtheorem{Def}[Prop]{Definition}

\newtheorem{As}[Prop]{Assumption}
\newtheorem{Test}[Prop]{Test}

\newtheorem{Cor}[Prop]{Corollary}

\newtheorem{Ex}[Prop]{Example}

\newtheorem*{Def*}{Definition}
\numberwithin{equation}{section}

\newcommand{\inas}{\operatorname{\stackrel{a.s.}{\to}}}
\newcommand{\inPrb}{\operatorname{\stackrel{\mathbb P}{\to}}}

\newcommand{\inD}{\operatorname{\stackrel{\mathcal{D}}{\to}}}

\newcommand{\Hess}{\mathrm{Hess}}

\newcommand{\Prb}{\mathbb{P}}

\newcommand{\EE}{\mathbb{E}}

\DeclareMathOperator*{\Cut}{Cut}

\DeclareMathOperator{\supp}{\mathrm{supp}}

\DeclareMathOperator{\cov}{\mathbf{cov}}

\title{\textbf{Moving-Basepoint Central Limit Theorems and
Corrected Wald Tests on Manifolds}}

\author[1]{Susovan Pal\footnote{susovan.pal@vub.be, susovan97@gmail.com}}
\affil[1]{\normalsize Mathematics and Data Science Group, Vrije Universiteit Brussel, Belgium}
\date{\today}

\begin{document}
\maketitle
\tableofcontents

\begin{abstract}
We investigate how CLTs on a Riemannian manifold and the resulting Wald-type $\chi^2$-tests are possibly affected by erroneously replacing a\textit{ moving} pooled sample Fréchet mean by a \textit{fixed}, population Fréchet mean while performing a Fréchet inference, a practice common in the relevant literature \citep{Huckemann2012, BhattacharyaBhattacharya2012}. In the case of arithmetic contrasts, we establish a new CLT  conforming to a computationally tractable moving chart and develop the corresponding Wald-type statistics. To prove the moving basepoint CLTs, we first generalize a CLT for general M-estimators from the i.i.d. setting of \citep{vanDerVaart1998} to finite grouped non-i.i.d. setting, and therefrom obtain a pooled Fréchet-mean CLT in the same setting, and then combine the latter with a suitable Taylor expansion of the parallel transport in \citep{HotzLeWood2024}. We show that the covariance in the CLT for the intrinsic contrasts does not change but the one in arithmetic contrast does when we work with moving charts. Under additional support conditions, we derive the corresponding, asymptotically $\chi^2$ Wald tests. A spherical counterexample shows that a previously claimed classical pooled-basepoint Hotelling limit can fail and illustrates the required correction. 
\end{abstract}

\section{Introduction}
\subsection{Motivation: inference in the tangent space at the pooled sample mean can differ from the same at a common population mean}

Two-sample Hotelling-type tests and one-way manifold MANOVA often compare
sample mean locations after representing them in a common tangent space. In manifold and shape analysis, a natural data-dependent choice of this space is the tangent space at a pooled sample mean. Examples
include the Procrustes tangent-space procedures of
\citet[Chapter~7]{DrydenMardia1998} and
\citet[Section~9.1.2]{DrydenMardia2016}, and the intrinsic manifold
MANOVA framework of \citet{HuckemannHotzMunk2010}. Procrustes tangent
coordinates are \textit{shape-specific} and should not, in general, be identified
with Riemannian logarithms. In the setting of Riemannian intrinsic means,
\citet[Section~6.1, Theorem~10]{Huckemann2012} study a two-sample
Hotelling $T^2$ statistic formed from individual tangent-space coordinates
at a pooled sample Fr\'echet mean. These procedures motivate studying
statistics whose common tangent-space basepoint is estimated from the
same observations used to form the statistic.

The corresponding \textit{asymptotic step} deserves attention. In the \textit{fixed-chart}
two-sample discussion of
\citet[Remark~2.6]{BhattacharyaPatrangenaru2005}, the chart does not move
with the sample. In contrast,
\citet[Section~5.4.1, equations~(5.19)--(5.21)]{BhattacharyaBhattacharya2012}
first obtain a two-sample intrinsic CLT in the normal chart at the \textit{unknown} intrinsic 
population Fréchet mean, and then\textit{ replace }this chart by one centered at the \textit{known}
pooled sample intrinsic Fréchet mean, owing to \textit{consistenc}y of that estimator.
But consistency alone does not, for a general statistic, determine whether
this replacement is negligible at the \textit{central-limit scale}. Likewise, the
proof of \citet[Theorem~10]{Huckemann2012} invokes an
$O_p(n^{-1/2})$ comparison between individual fixed- and pooled-basepoint coordinates. An estimate of this order does not, by itself, imply $o_p(1)$ after multiplication by $\sqrt n, n$ denoting the sample size; additional cancellation or an analysis of the full test statistic may be needed. This does not imply that all existing moving-basepoint procedures are invalid. However, for intrinsic logarithmic coordinates,
Subsection~\ref{subsec:huckemann-hotelling} gives a spherical
counterexample to the conventional chi-squared interpretation of
\citet[Theorem~10]{Huckemann2012} under its balanced-sample or covariance-equality
condition. The correction in Theorem~\ref{thm:moving-clts}\textup{(ii)}
explains the discrepancy, and Corollary~\ref{cor:wald} supplies
a correctly studentized test.

\subsection{The moving-basepoint problem, old vs. new CLT, implementation}

We consider $k$ mutually independent groups of manifold-valued data, with $X_1^{(r)},\ldots,X_{n_r}^{(r)}$ i.i.d. with law $P_r$ in group $r$,
$N:=\sum_{r=1}^k n_r$, and $n_r/N\to\lambda_r\in(0,1)$.
Under the null hypothesis $H_0$, the $P_r$ share a unique intrinsic population Fr\'echet mean $\mu$, though the distributions themselves
may differ. Write $X^{(r)}\sim P_r$. Let $\widehat\mu_{r,n_r}$ be the groupwise sample Fr\'echet means and $\widehat\mu_N$ the pooled sample Fr\'echet mean, both assumed to be unique. An \textit{implementable} intrinsic contrast between groups $r$ and $k$ is
$$
\log_{\widehat\mu_N}\widehat\mu_{r,n_r}
-\log_{\widehat\mu_N}\widehat\mu_{k,n_k},
\qquad r=1,\ldots,k-1.
$$
This vector belongs to the tangent space $T_{\widehat\mu_N}M$ varying with the total sample size $N$, and not to the \textit{fixed space}
$T_\mu M$ in which standard intrinsic Fr\'echet-mean CLTs are formulated (see \citet{BhattacharyaPatrangenaru2005} and \citet{KendallLe2011}).
Moreover, $T_{\widehat\mu_N}M$ varies with $N$. We therefore apply parallel transport, which is a linear isometry between tangent spaces preserving inner products and orthonormal bases, to the contrast along the \textit{unique }short minimizing geodesic from $\widehat\mu_N$ to $\mu$, using
$\tau_{\widehat\mu_N,\mu}:T_{\widehat\mu_N}M\to T_\mu M$.
This places the transported contrasts in a single fixed tangent space, where a joint Gaussian limit can be stated and compared with the fixed-basepoint limit, see Theorem \ref{thm:moving-clts}\textup{(i)}.

This puts both types of transported contrasts in the fixed tangent space
$T_\mu M$, where their limiting distributions can be compared with
those at the population basepoint. For clarity, consider the two-group
case $n_1/N\to\lambda\in(0,1)$ and define the \textit{arithmetic contrast}
\(
A_N(p):=\frac1{n_1}\sum_{j=1}^{n_1}\log_p X_j^{(1)}
-\frac1{n_2}\sum_{j=1}^{n_2}\log_p X_j^{(2)}.
\)
Under $H_0$ and the hypotheses of Theorem~\ref{thm:moving-clts}\textup{(ii)} for two groups ($k:=2$), we have two different CLT's, the first one being old and the second one new:
\[
\sqrt N\,A_N(\mu)
\xrightarrow{\mathcal D}N_{T_\mu M}(0,\Gamma_{\mathrm{fix}}),
\qquad
\sqrt N\,\tau_{\widehat\mu_N,\mu}A_N(\widehat\mu_N)
\xrightarrow{\mathcal D}N_{T_\mu M}(0,\Gamma^{(a)}),
\]
where the Fréchet function for the $r$-th group is $F_r(p):=\mathbb E[d_g(p,X^{(r)})^2]$, the tangential covariance is
$\Sigma_r:=\operatorname{Cov}(\log_\mu X^{(r)})$, the $(1,1)$-Hessian of $F_r$ is
$H_r:=\operatorname{Hess}_\mu^\sharp F_r$ and the affine combination of the groupwise Fréchet Hessian is 
$H_\lambda:=\lambda H_1+(1-\lambda)H_2$, and finally the old and the new covariances in the old and the new CLT for arithmetic contrasts between two groups ($k:=2$) are respectively:
\[
\Gamma_{\mathrm{fix}}
=\frac{\Sigma_1}{\lambda}+\frac{\Sigma_2}{1-\lambda},
\qquad
\Gamma^{(a)}
=\frac{H_2H_\lambda^{-1}\Sigma_1H_\lambda^{-1}H_2}{\lambda}
+\frac{H_1H_\lambda^{-1}\Sigma_2H_\lambda^{-1}H_1}{1-\lambda}.
\]
These covariance operators generally differ, even under $H_0$, although they coincide when $H_1=H_2$. \textit{Thus, using the fixed-basepoint covariance for the moving-basepoint arithmetic contrast need not yield a correctly calibrated Wald statistic.} In contrast, the intrinsic contrast of groupwise sample Fr\'echet means retains \emph{the same}
Gaussian limiting covariance at the fixed and moving basepoints after parallel transport (Theorem~\ref{thm:moving-clts}\textup{(i)}).
The corresponding $k$-group statements are given in
Theorem~\ref{thm:moving-clts}.

\subsubsection{The new arithmetic contrast CLT uses the\textit{ unknown} \texorpdfstring{$\mu$}{mu}, but the corresponding new Wald statistic does not} 
For implementation, however, the statistics and their plug-in covariance operators are represented \textit{directly} in $T_{\widehat\mu_N}M$ using an orthonormal basis of the same, and the Wald or Hotelling-type statistic, which is the inner product used in the tests are used in the moving tangent space {$T_{\widehat\mu_N}M$} itself, see Corollary \ref{cor:wald}. Their computation therefore does not require knowledge of the population Fr\'echet mean $\mu$.

\subsubsection{The issue of using one \textit{known} chart, possibly other than \texorpdfstring{$\operatorname{log_{\mu}}$}{log(mu)} since \texorpdfstring{$\mu$}{mu} is unknown} One might instead use a \textit{single fixed smooth chart }$\phi$ defined on a neighborhood of $\mu$, without the knowledge of $\mu$. This is a valid
alternative for comparing groupwise sample Fr\'echet
means. Indeed, the groupwise CLTs
\eqref{eq:groupwise-clt}, together with the multivariate
delta method applied to $\phi$, yield a joint Gaussian
limit for their \textit{fixed-chart contrasts}. However, this
does not establish the limiting distribution of the
different statistic obtained using the \textit{random chart}
$\log_{\widehat\mu_N}$. Nor can arithmetic averages of individual logarithms generally be replaced by averages of arbitrary
fixed-chart coordinates. Even when the coordinates
and their expectations are well defined, equality
of \textit{intrinsic Fr\'echet means} does not imply that in the chart $\phi,$ the Euclidean expectations are the same, i.e.
\(
\mathbb E[\phi(X^{(1)})]
=\cdots=
\mathbb E[\phi(X^{(k)})].
\)
{More explicitly, suppose that $\phi(X^{(r)})$ is defined $P_r$-almost surely and has a finite second moment for each group. Denoting by $Z_j^{(r)}:=\phi(X_j^{(r)})$, $\bar Z_r:=n_r^{-1}\sum_{j=1}^{n_r}Z_j^{(r)}$, $m_r^\phi:=\EE[Z_1^{(r)}]$, and $\Sigma_r^\phi:=\cov(Z_1^{(r)})$, the ordinary multivariate CLT gives 

$$\sqrt{n_r}(\bar Z_r-m_r^\phi)\xrightarrow{\mathcal D}\mathcal N_d(0,\Sigma_r^\phi).$$

By independence, for distinct groups $r,s$, the centered arithmetic contrast $\sqrt N[(\bar Z_r-\bar Z_s)-(m_r^\phi-m_s^\phi)]$ has a Gaussian limit with covariance $\Sigma_r^\phi/\lambda_r+\Sigma_s^\phi/\lambda_s$. Since $m_r^\phi-m_s^\phi$ need not vanish under $H_0:$ equality of the \textit{intrinsic} Fréchet means, the \textit{uncentered} contrast need not have a centered Gaussian limit, and the corresponding Hotelling statistic need not have its usual null distribution for testing equality of intrinsic Fr\'echet means. In contrast, the population logarithmic scores are centered under the hypotheses of Theorem~\ref{thm:moving-clts}\textup{(ii)}, while the arithmetic moving-logarithm contrast requires the Hessian-dependent moving-basepoint correction established therein.}
Consequently, such coordinate averages may test a
different population hypothesis.

\subsection{Contribution of the paper}

This paper produces CLTs for intrinsic and arithmetic contrasts with \textit{moving} basepoint as required for practical implementation without the knowledge of the common population Fréchet mean $\mu$ under $H_0$; \corrected{for intrinsic contrasts, }we do not impose any additional support or cut-locus restriction. On the way to prove the new CLTs in Theorem \ref{thm:moving-clts}; we needed to prove {a \textit{finite-group} M-estimator CLT (Proposition \ref{prop:grouped-m-estimator}) and its pooled-mean specialization (Proposition \ref{prop:pooled-clt}), with the linear representation needed \textit{when the group distributions differ}, under our population-level assumptions;} Even though these two propositions/CLT's are part of the proof of Theorem \ref{thm:moving-clts}, we stated them separately before that proof because we believe they could be of independent interest. We also derive an explicit first-order correction, the corresponding covariance operator, and a new moving-basepoint CLT for arithmetic tangent averages while justifying why the CLT does not change  for the intrinsic contrasts while replacing the fixed by moving basepoints: see Theorem \ref{thm:moving-clts} for both.
In contrast to the intrinsic result, the arithmetic CLT still requires either compact support away from the cut locus (Assumption \ref{as:compact-support}) or the compact-manifold cut-locus conditions of Lemma~\ref{lem:empirical-moving-log},
including the explicit empirical equicontinuity assumption.
Under Assumption~\ref{as:compact-support} (support away from the cut locus) and the further nondegeneracy conditions of the Wald result, we also obtain consistent plug-in covariances and asymptotic Wald tests; see
Lemma~\ref{lem:cov-consistency} and Corollary~\ref{cor:wald}.

The reason for imposing Assumption~\ref{as:compact-support} for plug-in estimation is that, without it, the population Hessian contains the cut-locus correction $-2J_{P_r}(\mu)$; see \citet[Theorem~1, equations~(14)--(15)]{HotzLeWood2024}. This correction is a surface integral over $\mathcal H_\mu$, an at most countable disjoint union of codimension-one hypersurfaces, involving their normal directions, the geometric factors $\kappa$ and $\tau'_{y|\mu}(0)$, and the density $\psi_r$; see \citet[Lemma~1 and equation~(20)]{HotzLeWood2024}. Identifying these hypersurfaces and consistently estimating the associated geometric and density terms is generally difficult, and no generally applicable consistent estimator is known to the author. Thus, the support restriction is imposed for the plug-in covariance estimation, not as a general requirement for the CLTs.

\subsection{Organization of the paper}

Section~\ref{sec:setup} introduces the multi-sample Riemannian
framework, notation, and assumptions.
Section~\ref{sec:pooled} proves strong consistency of the pooled and groupwise sample Fr\'echet means, develops a finite-group extension (Proposition \ref{prop:grouped-m-estimator}) of the classical $M$-estimator CLT \citep{vanDerVaart1998}, Theorem 5.23, and derives the pooled-mean CLT (Proposition \ref{prop:pooled-clt}) and the fixed-basepoint intrinsic contrast CLT, the last of which is \textit{not} new  at least when $k:=2$ but is included for completeness as we work with $k \ge 2$ groups. Section~\ref{sec:moving-clts} combines the intrinsic and arithmetic moving-basepoint CLTs. Subsection~\ref{sec:moving-log-linearization} states the empirical moving-logarithm linearization, with its proof in Appendix~\ref{app:empirical-moving-log}; Appendix~\ref{app:intrinsic-cancellation}
proves the stronger intrinsic cancellation estimate.
Finally, Section~\ref{sec:wald} establishes consistent plug-in
covariance estimation and the resulting asymptotic chi-squared Wald laws, with a corrected covariance for the arithmetic contrasts.
It then compares the arithmetic result with an earlier two-sample
Hotelling theorem, providing a spherical counterexample to its
classical limiting distribution.

\section{Setup}\label{sec:setup}
Let $M$ be a\textit{ complete} connected $d$-dimensional Riemannian manifold.  For
$r=1,\ldots,k$, let $X^{(r)}_1,X^{(r)}_2,\ldots$ be i.i.d. with law $P_r$, and
assume that the $k$ sequences are mutually independent.  We consider a
deterministic asymptotic regime of sample sizes $n_1,\ldots,n_k$ with
\[
 N:=\sum_{r=1}^k n_r\longrightarrow\infty,\qquad
 \lambda_{r,N}:=\frac{n_r}{N}\longrightarrow\lambda_r\in(0,1),\qquad
 \lambda:=(\lambda_1 \cdots \lambda_k), \qquad \sum_{r=1}^k\lambda_r=1.
\]
Write the population and sample Fréchet functions for the $r$-th group as
\[
 F_r(p):=\EE[d(p,X^{(r)})^2],\qquad
 F_{r,n_r}(p):=\frac1{n_r}\sum_{j=1}^{n_r}d(p,X^{(r)}_j)^2, \qquad
 1 \le r \le k.
\]
Under $H_0$, assume that the $F_r$ have the \textit{same unique} Fr\'echet mean $\mu$.
Moment, differentiability and nondegeneracy conditions are stated in the individual results below; compact support is not a standing assumption. Whenever a population quadratic expansion of $F_r$ at $\mu$ is assumed, $H_r$ denotes its Hessian at $\mu$, as in Proposition~\ref{prop:pooled-clt}.
When $F_r$ is twice differentiable at $\mu$, this means
$\langle H_ru,v\rangle_\mu=\Hess_\mu F_r(u,v)$.
Throughout the paper, $\log_a x$ denotes the usual minimizing Riemannian logarithm when $x\notin\Cut(a)$ and the zero vector in $T_aM$ when $x\in\Cut(a)$, following \citet[equation~(3)]{HotzLeWood2024}. This zero extension is a measurable convention, not an assertion of smoothness on the cut locus. On a normal neighborhood, $\log_a$ has its usual meaning. We refer explicitly to the \emph{ordinary, unextended logarithm} when no value on the cut locus is assigned. Likewise, $\Sigma_r:=\cov(\log_\mu X^{(r)})$ is used when the ordinary logarithm is defined $P_r$-almost surely and has a finite second moment. Put the affine Hessian and covariance:
$$
 H_\lambda:=\sum_{r=1}^k\lambda_rH_r,
 \qquad
 \Sigma_\lambda:=\sum_{r=1}^k\lambda_r\Sigma_r.
$$
Fix a relatively compact normal coordinate neighborhood $B$ of $\mu$ contained in $M\setminus\Cut(\mu)$. This serves only to define local coordinates and radial parallel transport; no geodesic convexity of $B$ or of the supports is assumed. For $a\in B$, let
$\tau_{a,\mu}:T_aM\to T_\mu M$ denote parallel transport
along the unique minimizing geodesic from $a$ to $\mu$.  We write $N_V(0,C)$ for the centered Gaussian law on a finite-dimensional inner-product space $V$ with covariance operator $C$.

\begin{As}[\textbf{Compact support away from the cut locus}]
\label{as:compact-support}
The set $C:=\bigcup_{r=1}^k\supp(P_r)$ is compact and $d_g(C,\Cut(\mu))>0$, with distance to the empty set interpreted as $+\infty$. No containment of $C$ in $B$ is assumed.
\end{As}
Assumption~\ref{as:compact-support} is imposed only where \textit{explicitly invoked}: as alternative \textup{(A)} of Lemma~\ref{lem:empirical-moving-log}
and for the plug-in covariance construction. The arithmetic
moving-basepoint result also permits alternative \textup{(B)} of that lemma. It is not required for pooled consistency,
the fixed-basepoint CLTs or the intrinsic moving-basepoint result.

\section{SLLN and CLT of the pooled sampled Fr\'echet mean}\label{sec:pooled}
Define
\[
 F_N:=\sum_{r=1}^k\lambda_{r,N}F_{r,n_r},\qquad
 F_\lambda:=\sum_{r=1}^k\lambda_rF_r,
\]
and let $\widehat\mu_N$ be a measurable global minimizer of $F_N$.
For each $r$, let $\widehat\mu_{r,n_r}$ be a measurable sequence of global minimizers of
$F_{r,n_r}$ that is the sample Fréchet mean for the group $r$. Denote by $\widehat\mu_N$ the global minimizer of $F_N$, which we call the sample Fréchet mean for the samples from all the groups combined. The goal of this section to establish a strong law of large numbers (SLLN) and a local tangent space central limit theorem (CLT) for the sample Fréchet mean $\widehat\mu_N$. Note that the CLT results of \citet{BhattacharyaPatrangenaru2005} do \textit{not} apply here to $\widehat{\mu}_N$ as our samples are \textit{merely independent, but not identical}.

\subsection{SLLN for the sample Fréchet mean \texorpdfstring{$\widehat\mu_N$}{mu-hat N} of all groups}
\begin{lem}[\textbf{Consistency of the pooled and groupwise sample Fr\'echet means}]
\label{lem:pooled-consistency}
Under the independent-group sampling assumptions, let $M$ be complete and
suppose that each $P_r$ has the same unique Fr\'echet mean $\mu$ and
\[
 \int_M d_g(\mu,x)^2\,dP_r(x)<\infty,
 \qquad r=1,\ldots,k.
\]
No compact-support or cut-locus assumption is required. Then
\[
 \widehat\mu_N\inas\mu,\qquad
 \widehat\mu_{r,n_r}\inas\mu\quad(r=1,\ldots,k).
\]
\end{lem}

\begin{proof}
Write
\[
 a_{r,n_r}:=\frac1{n_r}\sum_{j=1}^{n_r}d_g(\mu,X_j^{(r)}),
 \qquad a_r:=\int_M d_g(\mu,x)\,dP_r(x),
\]
and put $a_N:=\sum_{r=1}^k\lambda_{r,N}a_{r,n_r}$.
The moment assumption gives $a_r<\infty$. By the strong law of large
numbers, $a_{r,n_r}\to a_r$ almost surely, simultaneously for all
$r$, since $k$ is finite. In particular,
$a_N\to\sum_r\lambda_ra_r$ almost surely.

For $q,x\in M$ and $t=d_g(q,\mu)$, the reverse triangle inequality
implies
\[
 d_g(q,x)^2-d_g(\mu,x)^2\ge t^2-2t\,d_g(\mu,x).
\]
Thus
\[
 F_{r,n_r}(q)-F_{r,n_r}(\mu)\ge t^2-2t\,a_{r,n_r},
 \qquad
 F_N(q)-F_N(\mu)\ge t^2-2t\,a_N.
\]
Each empirical Fr\'echet function is continuous and, by these inequalities, takes values greater than its value at $\mu$ outside a sufficiently large closed ball. Since closed bounded balls are compact by Hopf--Rinow theorem for \textit{complete} manifolds (if incomplete, these balls need not be compact), each function attains a global minimum. Moreover, every groupwise minimizer satisfies
$d_g(\widehat\mu_{r,n_r},\mu)\le 2a_{r,n_r}$, and every pooled minimizer satisfies $d_g(\widehat\mu_N,\mu)\le2a_N$.
Choose a deterministic $R>2\max_{1\le r\le k}a_r+1$ and set
$K:=\overline B_g(\mu,R)$. Since
$a_N\le\max_r a_{r,n_r}$, almost surely all these minimizers lie
in the compact set $K$ for all sufficiently large $N$. For $q,z\in K$ and $x\in M$, the triangle inequality yields
\[
 |d_g(q,x)^2-d_g(z,x)^2|
 \le 2\bigl(d_g(\mu,x)+R\bigr)d_g(q,z).
\]
Consequently, $F_{r,n_r}$ and $F_r$ are Lipschitz on $K$ with
constants $2(a_{r,n_r}+R)$ and $2(a_r+R)$, respectively.
For any $\varepsilon>0$, compactness gives a finite
$\varepsilon$-net $z_1,\ldots,z_m$ of $K$, and hence
\begin{align*}
 \sup_{q\in K}|F_{r,n_r}(q)-F_r(q)|
 &\le \max_{1\le i\le m}|F_{r,n_r}(z_i)-F_r(z_i)|\\
 &\quad+2\varepsilon(a_{r,n_r}+a_r+2R).
\end{align*}
For each $z_i$, the squared distance is integrable because
$d_g(z_i,x)^2\le2d_g(\mu,x)^2+2R^2$, so the strong law makes the
maximum tend to zero almost surely. Taking
$\varepsilon=1/\ell$ for $\ell=1,2,\ldots$, and intersecting
the resulting countably many probability-one events, we obtain
from $a_{r,n_r}\to a_r$, simultaneously for all $r$,
\[
 \sup_{q\in K}|F_{r,n_r}(q)-F_r(q)|\longrightarrow0
 \qquad\text{almost surely}.
\]
Since $F_N=\sum_r\lambda_{r,N}F_{r,n_r}$ and
$F_\lambda=\sum_r\lambda_rF_r$, we have
\begin{align*}
 \sup_{q\in K}|F_N(q)-F_\lambda(q)|
 &\le \sum_{r=1}^k\lambda_{r,N}
  \sup_{q\in K}|F_{r,n_r}(q)-F_r(q)|\\
 &\quad+\sum_{r=1}^k|\lambda_{r,N}-\lambda_r|
  \sup_{q\in K}F_r(q)\longrightarrow0
  \qquad\text{almost surely},
\end{align*}
because $k$ is finite, $\lambda_{r,N}\to\lambda_r$,
and each $F_r$ is bounded on $K$.

Because each $F_r$ has the unique global minimizer $\mu$ and every
$\lambda_r>0$, $F_\lambda$ also has the unique global minimizer
$\mu$. For any $\delta>0$ for which
$K_\delta:=\{q\in K:d_g(q,\mu)\ge\delta\}$ is nonempty,
continuity and compactness give
\[
 \inf_{q\in K_\delta}
 \bigl(F_\lambda(q)-F_\lambda(\mu)\bigr)>0.
\]
The uniform convergence on $K$ and the fact that every pooled
minimizer eventually belongs to $K$ imply that no pooled minimizer
can belong to $K_\delta$ for all sufficiently large $N$, almost
surely. Thus $\widehat\mu_N\inas\mu$. Applying the same argument
to each $F_{r,n_r}$ and $F_r$ gives
$\widehat\mu_{r,n_r}\inas\mu$ for $r=1,\ldots,k$.
\end{proof}

In particular, with probability one, all sample means appearing below lie in
any fixed neighborhood of $\mu$ for all sufficiently large sample sizes.
This applies to $B$ and to the smaller neighborhoods used below.
Whenever a logarithm or a
parallel transport below is only canonically defined on that eventual event, we
may define the corresponding random vector arbitrarily on its complement; this
does not affect any asymptotic statement.


\subsection{CLT for the sample Fréchet mean \texorpdfstring{$\widehat\mu_N$}{mu-hat N} of all groups}

Our next goal is to establish a CLT for the sample Fréchet mean $\widehat{\mu}_N.$ We first extend the CLT for i.i.d. M-estimators
of \citet[Theorem~5.23, pp.~53--54]{vanDerVaart1998} to a \textit{fixed, finite }number
of \textit{independent} groups. Unlike a general triangular-array result, the
following proposition requires only groupwise regularity conditions.  In
particular, its criterion need not be a squared distance.

We extend \citet[Theorem~5.23, pp.~53--54]{vanDerVaart1998}
from i.i.d. observations to a fixed number of independent groups.
The criterion:cost/loss function below is general, as opposed to being convex or a function of the distance.

\begin{Prop}[\textbf{Finite-group extension of CLT for M-estimator \citet{vanDerVaart1998}, Theorem 5.23}]
\label{prop:grouped-m-estimator}
Let $\Theta\subset\mathbb R^d$ be open, and let
$\theta_0\in\Theta$. Fix $k\ge 1$. For $r=1,\ldots,k$, let
$X_1^{(r)},X_2^{(r)},\ldots X_{n_r}^{(r)}$ be i.i.d. with law $P_r$ on a common
measurable space $\mathcal X$, and assume that the $k$ sequences are
mutually independent. Let $n_r\to \infty$ for each group $r, N=\sum_{r=1}^k n_r\to\infty$ and
$\lambda_{r,N}:=n_r/N\to\lambda_r>0$ (so $\sum_r\lambda_r=1$).
For measurable criteria/cost/loss function $m_\theta:\mathcal X\to\mathbb R$, define:
\[
 \mathbb M_N(\theta):=\frac1N\sum_{r=1}^k\sum_{j=1}^{n_r}
          m_\theta(X_j^{(r)}),\quad
 M_r(\theta):=\int_{\mathcal X}m_\theta(x)\,dP_r(x),\quad
 M_N(\theta):=\sum_{r=1}^k\lambda_{r,N}M_r(\theta).
\]
Assume that:
\begin{enumerate}[label=\textup{(\roman*)},leftmargin=*]
\item For each $r$, $m_\theta$ is $P_r$-integrable for all
 $\theta\in\Theta$. In some neighborhood $U$ of $\theta_0$, on a
 $P_r$-full set independent of $\theta,\eta$, the map
 $\theta\mapsto m_\theta(x)$ is differentiable at $\theta_0$, with
 derivative $\psi(x):=\dot m_{\theta_0}(x)$, and
 \[
 |m_\theta(x)-m_\eta(x)|\le L_r(x)\|\theta-\eta\|
 \quad(\theta,\eta\in U),\qquad
 \int_{\mathcal X}L_r(x)^2\,dP_r(x)<\infty.
 \]
 So here $L_r$ is a nonnegative measurable, square-integrable envelope w.r.t. $P_r$.
\item Every groupwise population cost $M_r$ is maximized at $\theta_0$ and, for some symmetric
 matrix $V_r$, satisfies
 \[
 M_r(\theta_0+h)=M_r(\theta_0)
       +\tfrac12h^\top V_rh+o(\|h\|^2)\qquad(h\to0).
 \]
 The matrix $V_\lambda:=\sum_{r=1}^k\lambda_rV_r$ is nonsingular.
\item A measurable estimator $\widehat\theta_N\to\theta_0$ in
 probability is an approximate maximizer: there exist measurable
 $e_N\ge0$, with $e_N=o_p(N^{-1})$, such that
 \[
 \mathbb M_N(\widehat\theta_N)
 \ge\mathbb M_N(\theta)-e_N
 \qquad(\theta\in\Theta).
 \]
 This is the usual near-supremum condition of M-estimation when the
 supremum is measurable.
\end{enumerate}
Then $\int_{\mathcal X}\psi(x)\,dP_r(x)=0$ for every $r$, and, putting
$\Omega_\lambda:=\sum_{r=1}^k\lambda_r
 \int_{\mathcal X}\psi(x)\psi(x)^\top\,dP_r(x)$,
\[
 \sqrt N(\widehat\theta_N-\theta_0)
 =-V_\lambda^{-1}\frac1{\sqrt N}
   \sum_{r=1}^k\sum_{j=1}^{n_r}\psi(X_j^{(r)})+o_p(1),
 \qquad
 \sqrt N(\widehat\theta_N-\theta_0)
 \xrightarrow{\mathcal D}
 \mathcal N_d(0,V_\lambda^{-1}\Omega_\lambda V_\lambda^{-1}).
\]
\end{Prop}

\begin{proof}
We will use the standard notations from empirical process theory used in \cite{vanDerVaart1998}. For $P_r$-integrable $f$ put
\[
 \mathbb P_{r,n_r}f:=\frac1{n_r}\sum_{j=1}^{n_r}f(X_j^{(r)}),
 \quad P_rf:=\int_{\mathcal X}f(x)\,dP_r(x),\quad
 \mathbb G_{r,n_r}f:=\sqrt{n_r}(\mathbb P_{r,n_r}-P_r)f.
\]
Note that the last quantity $\mathbb G_{r,n_r}$ denotes the empirical process for the $r$-th group.
The local suprema used below have measurable versions: on the
full-measure Lipschitz sets in \textup{(i)}, they equal suprema over
fixed countable dense subsets of the coordinate balls.
Also write $\mathbb P_N:=\sum_r\lambda_{r,N}\mathbb P_{r,n_r}$,
$\overline P_N:=\sum_r\lambda_{r,N}P_r$, and
$\mathbb G_N:=\sqrt N(\mathbb P_N-\overline P_N)$.
Thus $\mathbb M_N(\theta)=\mathbb P_Nm_\theta$,
$M_N(\theta)=\overline P_Nm_\theta$, and
\begin{equation}\label{eq:grouped-decomposition}
 \mathbb G_Nf=\sum_{r=1}^k\sqrt{\lambda_{r,N}}\,
                  \mathbb G_{r,n_r}f.
\end{equation}
By (i), $\|\psi\|\le L_r$ $P_r$-a.s.; differentiation under the
integral and (ii) give $P_r\psi=\nabla M_r(\theta_0)=0$.
Each $V_r$ is negative semidefinite, hence their nonsingular positive
combination $V_\lambda$ is negative definite.

We follow the three steps of the proof of
\citet[Theorem~5.23, p.~54]{vanDerVaart1998}.
First, for every random vector $h_N=O_p(1)$, the grouped counterpart
of the stochastic approximation in that proof is
\begin{equation}\label{eq:grouped-stochastic-approx}
 \sqrt N\,\mathbb G_N
 (m_{\theta_0+h_N/\sqrt N}-m_{\theta_0})
 =h_N^\top\mathbb G_N\psi+o_p(1).
\end{equation}
Indeed, fix a group $r$ and set
$h_{r,N}:=\sqrt{\lambda_{r,N}}\,h_N=O_p(1)$.
Since $\theta_0+h_{r,N}/\sqrt{n_r}
=\theta_0+h_N/\sqrt N$, the i.i.d. result
\citet[Lemma~19.31, p.~284]{vanDerVaart1998}, applied to that group,
gives
\[
 \mathbb G_{r,n_r}\!\left[
 \sqrt{n_r}(m_{\theta_0+h_N/\sqrt N}-m_{\theta_0})
 -h_{r,N}^\top\psi\right]=o_p(1).
\]
The lemma applies even when $h_{r,N}$ depends on observations from
other groups, since it holds for arbitrary $O_p(1)$ random sequences.
Sum these relations over the fixed number $k$ of groups and use
\eqref{eq:grouped-decomposition} to obtain
\eqref{eq:grouped-stochastic-approx}.

Second, we establish $\sqrt N(\widehat\theta_N-\theta_0)=O_p(1)$,
as in \citet[Corollary~5.53, p.~77]{vanDerVaart1998}.
The groupwise Lipschitz condition, together with
\citet[Example~19.7, p.~271 and Corollary~19.35, p.~288]{vanDerVaart1998},
implies that, for sufficiently small $\delta>0$,
\[
 \mathbb E\sup_{\|h\|\le\delta}
 \left|\mathbb G_{r,n_r}
 (m_{\theta_0+h}-m_{\theta_0})\right|
 \le C_r\delta\quad(r=1,\ldots,k).
\]
Using \eqref{eq:grouped-decomposition} and the fact that $k$ is fixed,
\begin{equation}\label{eq:grouped-modulus}
 \mathbb E\sup_{\|h\|\le\delta}
 \left|\mathbb G_N(m_{\theta_0+h}-m_{\theta_0})\right|
 \le C\delta.
\end{equation}
Moreover, because $M_N=\sum_r\lambda_{r,N}M_r$, the expansions
in (ii) and $V_{\lambda,N}:=\sum_r\lambda_{r,N}V_r\to V_\lambda$
show that, for some $c,\delta_0>0$ and all sufficiently large $N$,
\begin{equation}\label{eq:grouped-gap}
 M_N(\theta_0+h)-M_N(\theta_0)
 \le-c\|h\|^2\qquad(\|h\|\le\delta_0).
\end{equation}
These are precisely the two bounds used in the proof of
\citet[Theorem~5.52, pp.~75--76]{vanDerVaart1998}, with
$\alpha=2$, $\beta=1$, and $n$ replaced by $N$.
For completeness, apply its shell argument to
$2^jA/\sqrt N<\|\widehat\theta_N-\theta_0\|
\le2^{j+1}A/\sqrt N$, within $\delta_0$.
On $\{Ne_N\le cA^2/2\}$, near-maximality and
\eqref{eq:grouped-gap} require a centered empirical fluctuation
of at least $c2^{2j}A^2/(2N)$.
By Markov's inequality and \eqref{eq:grouped-modulus},
the probability of the $j$th shell is at most $C'/(2^jA)$.
Summing over $j\ge0$, and using consistency and $Ne_N=o_p(1)$,
gives
\[
 \limsup_{N\to\infty}
 \mathbb P\!\left(\sqrt N\|\widehat\theta_N-\theta_0\|>A\right)
 \le\frac{C''}{A}.
\]
Hence $\sqrt N(\widehat\theta_N-\theta_0)=O_p(1)$.
This is the groupwise adaptation of Corollary~5.53 of \citet{vanDerVaart1998};
we have not applied that i.i.d. corollary directly to the pooled data.

Finally, put
\[
 h_N:=\sqrt N(\widehat\theta_N-\theta_0),\qquad
 S_N:=\mathbb G_N\psi
 =\frac1{\sqrt N}\sum_{r=1}^k\sum_{j=1}^{n_r}\psi(X_j^{(r)}).
\]
The independent groupwise central limit theorems imply
$S_N\xrightarrow{\mathcal D}\mathcal N_d(0,\Omega_\lambda)$;
in particular $S_N=O_p(1)$.
By (ii), the fixed number of population Taylor expansions, and
\eqref{eq:grouped-stochastic-approx}, for every random $t_N=O_p(1)$,
\begin{equation}\label{eq:grouped-local-quad}
 N\{\mathbb M_N(\theta_0+t_N/\sqrt N)
       -\mathbb M_N(\theta_0)\}
 =\tfrac12t_N^\top V_\lambda t_N+t_N^\top S_N+o_p(1).
\end{equation}
Equation~\eqref{eq:grouped-local-quad} holds for $t_N=h_N$ and for
$t_N^*:=-V_\lambda^{-1}S_N=O_p(1)$.
Because $\Theta$ is open, $\theta_0+t_N^*/\sqrt N\in\Theta$
with probability tending to one. Near-maximality and $Ne_N=o_p(1)$
therefore yield
\[
 \tfrac12h_N^\top V_\lambda h_N+h_N^\top S_N
 \ge \tfrac12(t_N^*)^\top V_\lambda t_N^*
       +(t_N^*)^\top S_N+o_p(1).
\]
Completing the square, just as in the proof of Theorem~5.23 in \citet{vanDerVaart1998},
\[
 \tfrac12(h_N-t_N^*)^\top V_\lambda(h_N-t_N^*)
 \ge-o_p(1).
\]
Since $V_\lambda$ is negative definite,
$h_N-t_N^*=o_p(1)$, proving the stated linear representation.
The Gaussian limit follows from the groupwise CLT and Slutsky's theorem.
\end{proof}

\subsubsection{Translation table from general \texorpdfstring{$M$}{M}-estimation to geometric statistics}
The following table specializes the general M-estimation quantities to
the geometric setting.  It uses normal coordinates at the common
population Fr\'echet mean $\mu$; the two columns are identities, not
additional assumptions. The symbols on the left column are common in the M-estimation literature, especially in \citet{vanDerVaart1998}, and the ones in the right columns are common in geometric statistic literature.
\begin{center}
\small
\renewcommand{\arraystretch}{1.18}
\begin{tabular}{@{}p{0.37\linewidth}p{0.56\linewidth}@{}}
\hline
\textbf{General grouped M-estimation} & \textbf{Pooled intrinsic Fr\'echet mean}\\
\hline
$\mathcal X$, $\theta_0$ & $M$, $0\in T_\mu M\cong\mathbb R^d$\\
$m_\theta(x)$ & $-d_g(\exp_\mu\theta,x)^2$\\
$M_r(\theta)$ & $-F_r(\exp_\mu\theta)$\\
$\mathbb M_N(\theta)$ & $-F_N(\exp_\mu\theta)$\\
$M_N(\theta)$ & $-\sum_{r=1}^k\lambda_{r,N}F_r(\exp_\mu\theta)$\\
$\widehat\theta_N$ & $\xi_N=\log_\mu\widehat\mu_N$\\
$\psi(x)$ & $2\log_\mu x$\\
$V_r$, $V_\lambda$ & $-H_r$, $-H_\lambda$\\
$\Omega_\lambda$ & $4\Sigma_\lambda$\\
\hline
\end{tabular}
\end{center}
In particular, $\mathbb M_N(\theta)=-F_N(\exp_\mu\theta)$,
whereas $M_N(\theta)$ is its expectation and need not equal the limiting
$-F_\lambda(\exp_\mu\theta)$.
For comparison, consistency in the classical i.i.d. setting is treated
by \citet[Theorem~5.7, p.~45]{vanDerVaart1998}; in our grouped setting
it was proved directly in Lemma~\ref{lem:pooled-consistency}.

For the normal-chart application below, fix a sufficiently small coordinate
ball $\Theta=B(0,\rho)\subset T_\mu M$ and set
$E_N:=\{\widehat\mu_N\in\exp_\mu\Theta\}$.
Define $\widehat\theta_N:=\log_\mu\widehat\mu_N$ on $E_N$ and
$\widehat\theta_N:=0$ otherwise. Lemma~\ref{lem:pooled-consistency}
gives $\Prb(E_N)\to1$ and $\widehat\theta_N\inPrb0$.
For $\mathbb M_N(\theta):=-F_N(\exp_\mu\theta)$, the measurable error
$e_N:=(\sup_{\theta\in\Theta}\mathbb M_N(\theta)
-\mathbb M_N(\widehat\theta_N))_+$ is finite and is zero on $E_N$,
because $\widehat\mu_N$ is a global minimizer. Hence
$\Prb(Ne_N>\varepsilon)\le\Prb(E_N^c)\to0$ for every $\varepsilon>0$.
These choices supply condition \textup{(iii)} of
Proposition~\ref{prop:grouped-m-estimator} without any support restriction.

\begin{Prop}[\textbf{Pooled-mean CLT and linearization}]
\label{prop:pooled-clt}
Suppose the independent groups have the same unique population
Fr\'echet mean $\mu$ and, for each $r=1,\ldots,k$,
\[
 \int_M d_g(\mu,x)^2\,dP_r(x)<\infty.
\]
Note that here $p\mapsto d_g(p,x)^2$ is differentiable at $p=\mu$
for $P_r$-almost every $x$ thanks to \citet[Theorem~1.4]{LytchakHuckemann2025}. Write $\log_\mu x$ for the initial
velocity of the unique minimizing geodesic from $\mu$ to $x$ on
this full-measure set, retaining the zero-extension convention on $\Cut(\mu)$.
Assume also that, as $h\to0$ in $T_\mu M$,
\(
 F_r(\exp_\mu h)
 =F_r(\mu)+\tfrac12\langle H_rh,h\rangle_\mu+o(\|h\|^2)
\)
for some self-adjoint operator $H_r:T_\mu M\to T_\mu M$.
Define
\[
 H_\lambda:=\sum_{r=1}^k\lambda_rH_r,\qquad
 \Sigma_r:=\operatorname{Cov}(\log_\mu X^{(r)}),\qquad
 \Sigma_\lambda:=\sum_{r=1}^k\lambda_r\Sigma_r,
\]
and assume that $H_\lambda$ is positive definite. No
individual $H_r$ is assumed nonsingular. No restriction on the
supports of the $P_r$ is imposed. Let $\widehat\mu_N$ be a measurable global minimizer of $F_N$ and put
\(
 U_{r,N}:=\frac1{n_r}\sum_{j=1}^{n_r}\log_\mu X^{(r)}_j,
 \, \xi_N:=\log_\mu\widehat\mu_N.
\)
The logarithm defining $\xi_N$ is understood in a fixed normal
neighborhood of $\mu$; set $\xi_N=0$ when $\widehat\mu_N$
is outside that neighborhood. Then
 \begin{equation}
 \sqrt N\,\xi_N
 =2H_\lambda^{-1}\sum_{r=1}^k\lambda_{r,N}\sqrt N\,U_{r,N}
   +o_p(1), \label{eq:pooled-linearization}
\end{equation}
and consequently
\begin{equation}
 \sqrt N\,\log_\mu\widehat\mu_N
 \xrightarrow{\mathcal D}N_{T_\mu M}\!\left(0,
  4H_\lambda^{-1}\Sigma_\lambda H_\lambda^{-1}\right). \label{eq:pooled-clt}
\end{equation}
\end{Prop}

\begin{proof}
Fix an orthonormal basis of $T_\mu M$. Choose an open ball
$\Theta=B(0,\rho)\subset T_\mu M\cong\mathbb R^d$ small enough
that $\exp_\mu$ is a diffeomorphism on a neighborhood of
$\overline\Theta$. For $\theta\in\Theta$ and $x\in M$, set
\[
 m_\theta(x):=-d_g(\exp_\mu\theta,x)^2,\qquad
 M_r(\theta):=\int_Mm_\theta(x)\,dP_r(x)
              =-F_r(\exp_\mu\theta).
\]
We check conditions \textup{(i)}--\textup{(iii)} of the preceding
finite-group M-estimation proposition.

\emph{Condition \textup{(i)}.}
For $\theta\in\Theta$, the triangle inequality gives
\[
 d_g(\exp_\mu\theta,x)^2
 \le 2d_g(\mu,x)^2+2\rho^2,
\]
so $m_\theta$ is integrable under every $P_r$. Since
$\exp_\mu$ is smooth on a neighborhood of $\overline\Theta$,
there is $K<\infty$ such that
\[
 d_g(\exp_\mu\theta,\exp_\mu\eta)
 \le K\|\theta-\eta\|\qquad(\theta,\eta\in\Theta).
\]
The triangle inequality then yields, for every $x\in M$,
\begin{align*}
 |m_\theta(x)-m_\eta(x)|
 &\le d_g(\exp_\mu\theta,\exp_\mu\eta)
       \bigl[d_g(\exp_\mu\theta,x)+d_g(\exp_\mu\eta,x)\bigr]\\
 &\le 2K\bigl(d_g(\mu,x)+\rho\bigr)\|\theta-\eta\|.
\end{align*}
The function $L(x):=2K(d_g(\mu,x)+\rho)$ satisfies
$\int_M L(x)^2\,dP_r(x)<\infty$ for every $r$.
This bound holds for all $\theta,\eta\in\Theta$ and every $x$;
it does not require differentiability at points other than $\mu$.
By the assumed differentiability and the first-variation formula,
\[
 \psi(x):=D_\theta m_\theta(x)\big|_{\theta=0}
          =2\log_\mu x\qquad(P_r\text{-a.s.}).
\]
Because the difference quotients are bounded by $L$, we may
pass their limits through the expectation. Hence
\[
 \int_M\psi(x)\,dP_r(x)=\nabla M_r(0)=0,
 \qquad
 \int_M\psi(x)\psi(x)^\top\,dP_r(x)=4\Sigma_r.
\]
The first equality uses the fact that $M_r$ is maximized at $0$.

\emph{Condition \textup{(ii)}.}
Each $M_r$ is maximized at $0$, and the assumed population expansion
implies
\[
 M_r(h)=M_r(0)+\tfrac12h^\top V_rh+o(\|h\|^2),
 \qquad V_r=-H_r.
\]
Thus $V_\lambda=\sum_r\lambda_rV_r=-H_\lambda$ is negative
definite. The operators $H_r$ are the coefficients of the
\emph{population} quadratic expansions; no interchange of a
second derivative and an expectation is used.

\emph{Condition \textup{(iii)}}: holds because of Lemma \ref{lem:pooled-consistency}.

Proposition~\ref{prop:grouped-m-estimator}, applied with $\theta_0=0$, $V_\lambda=-H_\lambda$
and $\Omega_\lambda=4\Sigma_\lambda$, now gives
\[
 \sqrt N\,\widehat\theta_N
 =H_\lambda^{-1}\frac1{\sqrt N}
       \sum_{r=1}^k\sum_{j=1}^{n_r}2\log_\mu X_j^{(r)}+o_p(1).
\]
Since $\widehat\theta_N=\xi_N$ on $E_N$, and since
\(
 \frac1{\sqrt N}\sum_{j=1}^{n_r}\log_\mu X_j^{(r)}
 =\lambda_{r,N}\sqrt N\,U_{r,N},
\)
we obtain \textup{\eqref{eq:pooled-linearization}}. The Gaussian limit in
Proposition~\ref{prop:grouped-m-estimator} gives \textup{\eqref{eq:pooled-clt}}.
\end{proof}

\begin{Remark}[\textbf{Comparison with Kendall--Le (2011), Theorem 4 (CLT for non-iid sequences of random variables)}]
\label{rem:kendall-le-comparison}
Theorem~4 of \citet[pp.~345--346]{KendallLe2011} does not
directly imply Proposition~\ref{prop:pooled-clt}.
In particular, Conditions~2 and~3 of their theorem require,
respectively, normalized control of the local variation and
second moments of the individual squared-distance Hessians
\[
 H_i(p):=\operatorname{Hess}_p
 \left(\tfrac12 d_g(p,X_i)^2\right);
\]
see their equation~(5.4), p.~343.
Neither condition is imposed in Proposition~\ref{prop:pooled-clt}.

Furthermore, their Condition~4 requires asymptotic
nonsingularity of
\[
 \widetilde H_n
 :=\frac{1}{2\phi_n(\mu)}
   \sum_{i=1}^n\mathbb E[H_i(\mu)],
\]
whereas Proposition~\ref{prop:pooled-clt} assumes positive
definiteness of
\(
 H_\lambda=\sum_{r=1}^k\lambda_r H_r,
\)
where each $H_r$ is determined by the second-order expansion
of the population Fr\'echet function $F_r$ at $\mu$.
These two types of Hessians need not coincide, since
interchanging the Hessian and expectation is not always
justified, particularly in the presence of cut-locus effects.
Neither result requires the individual Hessians to be
nonsingular.

Finally, provided that
\(
 v_\lambda:=\sum_{r=1}^k\lambda_rF_r(\mu)>0,
\)
the Lindeberg-type condition (Condition~5 of
\citet[Theorem~4]{KendallLe2011}) holds automatically
in our finite-group setting. Indeed, the aggregate
population variation in that paper, denoted by 
\(
 \phi_N(\mu):=\frac12\sum_{r=1}^k n_rF_r(\mu)
\)
satisfies $\phi_N(\mu)/N\to v_\lambda/2>0$.
Since the number of distribution classes is fixed
and each distribution has a finite second moment,
the required normalized tail expectations converge
to zero by dominated convergence.
If $v_\lambda=0$, every $P_r$ is concentrated at $\mu$;
our pooled-mean result then holds trivially with zero
limiting covariance, but Kendall--Le's positive-variation
assumption is not satisfied.

\end{Remark}

\subsubsection{Intrinsic and arithmetic contrasts and their \textit{fixed-basepoint} CLTs}

\paragraph{CLT for fixed-basepoint intrinsic contrasts:}
For $r<k$, define the $r$-th groupwise contrasts $\Delta_{r,N}^{(i)}:=\log_{\widehat\mu_N}\widehat\mu_{r,n_r}-\log_{\widehat\mu_N}\widehat\mu_{k,n_k}$, and collect these components as $\Delta_N^{(i)}:=(\Delta_{1,N}^{(i)},\ldots,\Delta_{k-1,N}^{(i)})\in(T_{\widehat\mu_N}M)^{k-1}$. The corresponding fixed-basepoint contrast is $\Delta_{\mu,N}^{(i)}:=(\log_\mu\widehat\mu_{r,n_r}-\log_\mu\widehat\mu_{k,n_k})_{r=1}^{k-1}$. Write $\mathcal T_N:=\tau_{\widehat\mu_N,\mu}^{\oplus(k-1)}$ for componentwise parallel transport to $(T_\mu M)^{k-1}$. These objects are defined on the event that all sample means lie in a common small normal neighborhood; define them arbitrarily on its complement, whose probability tends to zero by the SLLN in Lemma~\ref{lem:pooled-consistency}.

For the intrinsic contrasts, suppose that each group satisfies the moment, almost-sure differentiability and population expansion conditions of Proposition~\ref{prop:pooled-clt}, and that every $H_r$ is positive definite. Choose the groupwise estimators as measurable functions of their respective samples. Applying Proposition~\ref{prop:grouped-m-estimator} separately to each group, with the normal-coordinate specialization used in Proposition~\ref{prop:pooled-clt}, gives
\begin{equation}\label{eq:groupwise-clt}
\sqrt{n_r}\log_\mu\widehat\mu_{r,n_r}
\xrightarrow{\mathcal D}N_{T_\mu M}(0,Q_r),
\qquad Q_r:=4H_r^{-1}\Sigma_rH_r^{-1},\quad r=1,\ldots,k.
\end{equation}
The limits in \eqref{eq:groupwise-clt} are jointly independent, since the samples are independent. By the groupwise CLTs \eqref{eq:groupwise-clt}, independence and $n_r/N\to\lambda_r>0$ yield the following CLT for the fixed-basepoint intrinsic contrast:
\(
\sqrt N\,\Delta_{\mu,N}^{(i)}
\xrightarrow{\mathcal D}
(Y_1-Y_k,\ldots,Y_{k-1}-Y_k),
\)
where $Y_1,\ldots,Y_k$ are independent centered Gaussian vectors in $T_\mu M$ with $\operatorname{Cov}(Y_r)=Q_r/\lambda_r$, $r=1,\ldots,k$. Consequently, the limiting covariance operator $\Gamma^{(i)}:(T_\mu M)^{k-1}\to(T_\mu M)^{k-1}$ has the $(k-1)\times(k-1)$ block representation
\begin{equation}\label{eq:intrinsic-covariance}
\Gamma^{(i)}=\left(\Gamma_{rs}^{(i)}\right)_{1\le r,s<k},
\qquad
\Gamma_{rs}^{(i)}=\delta_{rs}\frac{Q_r}{\lambda_r}+\frac{Q_k}{\lambda_k}.
\end{equation}
Indeed, since $Y_1,\ldots,Y_k$ are independent, the diagonal blocks are $\Gamma_{rr}^{(i)}=\operatorname{Cov}(Y_r-Y_k,Y_r-Y_k)=Q_r/\lambda_r+Q_k/\lambda_k$, whereas, for $r\neq s$, the off-diagonal blocks are $\Gamma_{rs}^{(i)}=\operatorname{Cov}(Y_r-Y_k,Y_s-Y_k)=\operatorname{Cov}(Y_k,Y_k)=Q_k/\lambda_k$. Thus, the off-diagonal covariance arises because all contrasts share the same reference group $k$.

\paragraph{CLT for fixed-basepoint arithmetic contrasts:}
We next consider the fixed-basepoint arithmetic contrasts $(U_{r,N}-U_{k,N})_{r=1}^{k-1}$, where $U_{r,N}$ is defined in Proposition~\ref{prop:pooled-clt}. Under $H_0$, the first-variation argument in that proposition gives $\mathbb E[\log_\mu X^{(r)}]=0$, while the finite second-moment assumption gives $\operatorname{Cov}(\log_\mu X^{(r)})=\Sigma_r$. The ordinary multivariate CLT, independence of the groups, and $n_r/N\to\lambda_r>0$ therefore yield
\[
\sqrt N\,(U_{r,N}-U_{k,N})_{r=1}^{k-1}
\xrightarrow{\mathcal D}
N_{(T_\mu M)^{k-1}}\left(0,
\left(\delta_{rs}\frac{\Sigma_r}{\lambda_r}
+\frac{\Sigma_k}{\lambda_k}\right)_{1\le r,s<k}\right).
\]
The limiting covariance operator on $(T_\mu M)^{k-1}$ has the same $(k-1)\times(k-1)$ block structure as \eqref{eq:intrinsic-covariance}, with $Q_r$ replaced by $\Sigma_r$. In particular, the diagonal blocks are $\Sigma_r/\lambda_r+\Sigma_k/\lambda_k$, while the off-diagonal blocks are $\Sigma_k/\lambda_k$, since all contrasts share the same reference group $k$. This fixed-basepoint arithmetic CLT requires no positive-definiteness assumption on the individual Hessians $H_r$.


\section{Moving-basepoint central limit theorems}\label{sec:moving-clts}
Fix $k\ge2$ and retain the intrinsic contrasts $\Delta_N^{(i)}$, their fixed-basepoint counterparts $\Delta_{\mu,N}^{(i)}$, and componentwise parallel transport $\mathcal T_N$ from the preceding section. Use the logarithm convention of Section~\ref{sec:setup}. The arithmetic contrasts and their covariance operators are defined by
\begin{equation}\label{eq:arithmetic-averages}
\begin{aligned}
V^{(a)}_{r,N}
&:=\frac1{n_r}\sum_{j=1}^{n_r}
\log_{\widehat\mu_N}(X_j^{(r)}) \in T_{\widehat\mu_N}M,\\
\Delta_N^{(a)}
&:=\bigl(V^{(a)}_{1,N}-V^{(a)}_{k,N},
\ldots,V^{(a)}_{k-1,N}-V^{(a)}_{k,N}\bigr) \in (T_{\widehat\mu_N}M)^{k-1}.
\end{aligned}
\end{equation}
Under Assumption~\ref{as:compact-support}, these
coincide with the ordinary logarithm averages
with probability tending to one. For the moving-basepoint CLT for arithmetic contrasts in Theorem \ref{thm:moving-clts}(ii), when $H_\lambda$ is positive definite, put
\begin{equation}\label{eq:arithmetic-blocks}
\begin{aligned}
B_{rs}
&:=(\delta_{rs}-\delta_{ks})I_d
   -\lambda_s(H_r-H_k)H_\lambda^{-1},
\quad 1\le r \le k-1,\quad 1\le s\le k,\\
B&:=(B_{rs})_{r\le k-1,\,s\le k}: (T_\mu M)^k\to(T_\mu M)^{k-1},\\
D_\lambda&:=\operatorname{diag}\!\left(
\frac{\Sigma_1}{\lambda_1},\ldots,
\frac{\Sigma_k}{\lambda_k}\right):
(T_\mu M)^k\to(T_\mu M)^k.
\end{aligned}
\end{equation}
Here $I_d$ is the identity on $T_\mu M$,
$\delta_{rs}$ is the Kronecker delta, and
$B_{rs}:T_\mu M\to T_\mu M$ is a linear operator.
Thus $B:(T_\mu M)^k\to(T_\mu M)^{k-1}$
is a $(k-1)\times k$ block operator, with each block
of dimension $d\times d$, while
$D_\lambda:(T_\mu M)^k\to(T_\mu M)^k$
is block diagonal. The operator $B$ collects the
coefficients of the first-order linearization of the
arithmetic contrasts, including the moving-basepoint
correction, as established in \eqref{eq:arithmetic-linearization}.
Here $B^*$ denotes the adjoint of $B$.
Define the covariance operator $\Gamma^{(a)}$ on $(T_\mu M)^{k-1}$ by
\begin{equation}\label{eq:arithmetic-covariance}
\Gamma^{(a)}:=BD_\lambda B^*,
\qquad
\Gamma_{rt}^{(a)}
=\sum_{s=1}^k
B_{rs}\frac{\Sigma_s}{\lambda_s}B_{ts}^*,
\quad 1\le r,t<k.
\end{equation}

\subsection{Statements of the new, moving-basepoint CLT's}

{\textit{Local sampling regime.} Let the laws $P_1,\ldots,P_k$ in Section \ref{sec:setup} share the unique population Fr\'echet mean $\mu$. Fix bounded measurable score functions $s_r:M\to\mathbb R$ with $\int_M s_r\,dP_r=0$ of the one-dimensional submodels defined below. For all sufficiently large $N$, let the observations within group $r$ be i.i.d. with law $P_{r,N}$ defined by $dP_{r,N}/dP_r=1+s_r/\sqrt N$, independently across groups. These are probability laws: their densities $dP_{r,N}/dP_r$ integrate to one because $\int_M s_r\,dP_r=0$, and lie between $1/2$ and $3/2$ for sufficiently large $N$. All population quantities and regularity assumptions below refer to the baseline laws $P_1,\ldots,P_k$; the sample means and other statistics are the same measurable functions of the observations as before. Write $\Prb_{0,N},\EE_{0,N}$ and $\Prb_{s,N},\EE_{s,N}$ 
for the \textit{joint} probability laws and corresponding 
expectations of all $N$ observations under baseline and 
locally perturbed sampling induced by the score functions 
$s_1,\ldots,s_k$, respectively. By independence,
\(
\Prb_{0,N}=\bigotimes_{r=1}^{k}P_r^{\otimes n_r},
\,
\Prb_{s,N}=\bigotimes_{r=1}^{k}P_{r,N}^{\otimes n_r},
\)
where $\otimes$ denotes the product of probability measures. Under the applicable moment and differentiability hypotheses, define $a_r:=\int_M(\log_\mu x)s_r(x)\,dP_r(x)\in T_\mu M$. For baseline sampling take every $s_r\equiv0$, so $a_r=0$; this recovers the original null model under $H_0$ exactly.}

\begin{Def}[\textbf{Contiguous alternatives}]\label{def:contiguous alternatives}
The alternative sampling regime following the law $\Prb_{s,N}$ is \emph{contiguous} with
respect to the baseline regime $\Prb_{0,N}$ if, for every sequence of measurable sample events $A_N$, $\Prb_{0,N}(A_N)\to0$
implies $\Prb_{s,N}(A_N)\to0$. The two regimes are
\emph{mutually contiguous} if the reverse implication
also holds \citep[Definition~6.3, p.~87]{vanDerVaart1998}.
\end{Def}

The baseline and local sampling regimes introduced above are indeed mutually contiguous. To see this, let $L_N:=d\Prb_{s,N}/d\Prb_{0,N}=\prod_{r=1}^k\prod_{j=1}^{n_r}(1+s_r(X_j^{(r)})/\sqrt N)$ denote the joint likelihood ratio. By independence and boundedness of the centered scores, for all sufficiently large $N$,
\[
\EE_{0,N}[L_N^2]\le e^C,\qquad \EE_{s,N}[L_N^{-2}]\le e^{2C},
\]
where $C:=\max_{1\le r\le k}\int_M s_r^2\,dP_r$. Mutual contiguity follows by the Cauchy--Schwarz inequality. Consequently, consistency in probability, $O_p(1)$ bounds, and $o_p(1)$ remainders transfer in both directions for the same statistics and deterministic scalings.

\begingroup
\color{black}
\begin{Th}[\textbf{Moving-basepoint central limit theorems}]\label{thm:moving-clts}
\color{black}
Under $H_0$ (group $r$ sampled from $P_r$) or the local sampling regime above (group $r$ sampled from $P_{r,N}$) and the independent-group sampling regime $n_r/N\to\lambda_r\in(0,1)$, the following conclusions hold under their respective hypotheses.
\begin{enumerate}[label=\textup{(\roman*)},leftmargin=*]
\item \textbf{Intrinsic contrasts.} Under the groupwise conditions stated for \eqref{eq:groupwise-clt}, including positive definiteness of every $H_r$,
\begin{equation}\label{eq:intrinsic-moving-clt}
\sqrt N\,\mathcal T_N\Delta_N^{(i)}\xrightarrow{\mathcal D}N_{(T_\mu M)^{k-1}}({h^{(i)}},\Gamma^{(i)}).
\end{equation}

Here $h^{(i)}:=(2H_r^{-1}a_r-2H_k^{-1}a_k)_{r=1}^{k-1}$, and $\Gamma^{(i)}$ is a $(k-1)\times(k-1)$ block covariance operator with blocks
$\Gamma_{rs}^{(i)}=\delta_{rs}Q_r/\lambda_r+Q_k/\lambda_k$,
$1\le r,s<k$, where $Q_r:=4H_r^{-1}\Sigma_rH_r^{-1}$; see \eqref{eq:intrinsic-covariance}.
The basepoint replacement requires only pooled-mean consistency; neither Assumption~\ref{as:compact-support}, Lemma~\ref{lem:empirical-moving-log}, nor a pooled-mean CLT is required.

\item \textbf{Arithmetic contrasts.} Under the hypotheses of Proposition~\ref{prop:pooled-clt} and either alternative of Lemma~\ref{lem:empirical-moving-log},
\begin{equation}\label{eq:arithmetic-moving-clt}
\sqrt N\,\mathcal T_N\Delta_N^{(a)}\xrightarrow{\mathcal D}N_{(T_\mu M)^{k-1}}({h^{(a)}},\Gamma^{(a)}).
\end{equation}
Here $h^{(a)}:=B(a_1,\ldots,a_k)
=\left(\sum_{s=1}^{k}B_{rs}a_s\right)_{r=1}^{k-1}
\in(T_\mu M)^{k-1}$. Here $\Gamma^{(a)}=BD_\lambda B^*$ is given by \eqref{eq:arithmetic-blocks}--\eqref{eq:arithmetic-covariance}, using the full population Hessians $H_r$. Positive definiteness is required only of $H_\lambda$, not of individual $H_r$ for any $r$. 
\end{enumerate}
\end{Th}
\endgroup

\begin{Cor}[\textbf{Case $k=2$}]\label{cor:two-group-moving-clts}
Under $H_0$ or the local sampling regime of Theorem~\ref{thm:moving-clts}, let $n_1:=n$, $n_2:=m$, $N=n+m$, and $n/N\to\lambda\in(0,1)$. Under the hypotheses of the Theorem~\ref{thm:moving-clts} above, we have the two following moving basepoint CLT's:

\begin{equation}\label{eq:moving-basept-CLT-2-groups}
    \begin{aligned}
\textbf{intrinsic:}
\quad \sqrt N\,\tau_{\widehat\mu_N,\mu}\Delta_N^{(i)}&\inD N_{T_\mu M}({h^{(i)}},\Gamma^{(i)}),\\
\textbf{arithmetic:}
\quad\sqrt N\,\tau_{\widehat\mu_N,\mu}(V_{1,N}^{(a)}-V_{2,N}^{(a)})&\inD N_{T_\mu M}({h^{(a)}},\Gamma^{(a)}),
   \end{aligned}
\end{equation}
where the \textit{intrinsic} and \textit{arithmetic limiting covariances} are as below:
$\Gamma^{(i)}=4H_1^{-1}\Sigma_1H_1^{-1}/\lambda+4H_2^{-1}\Sigma_2H_2^{-1}/(1-\lambda)$ and $\Gamma^{(a)}=H_2H_\lambda^{-1}\Sigma_1H_\lambda^{-1}H_2/\lambda+H_1H_\lambda^{-1}\Sigma_2H_\lambda^{-1}H_1/(1-\lambda)$, with $H_\lambda=\lambda H_1+(1-\lambda)H_2$. The means are $h^{(i)}=2H_1^{-1}a_1-2H_2^{-1}a_2$ and $h^{(a)}=H_2H_\lambda^{-1}a_1-H_1H_\lambda^{-1}a_2$; both vanish when every $s_r\equiv0$. These formulas follow from Theorem~\ref{thm:moving-clts} and the two-group identities for $B$ in Remark~\ref{rem:two-groups}.
\end{Cor}

\subsection{Proof of the moving basepoint CLT's in Theorem \ref{thm:moving-clts}}\label{subsec:proof-moving-clts}
The proof of part (i) of Theorem \ref{thm:moving-clts} is in the subsection \ref{subsec:Common Hessian expansion and intrinsic cancellation}. The proof of part (ii) is in the subsection \ref{sec:moving-log-linearization}, which involves the technical Lemma \ref{lem:empirical-moving-log}, that is proved in the Appendix \ref{app:empirical-moving-log}.

Before proceeding to the two proofs, we establish two
probabilistic results common to both: the \textit{mutual
contiguity} (cf. Definition \ref{def:contiguous alternatives}) of the baseline and locally perturbed sampling
laws, which permits the \textit{transfer of stochastic remainder
estimates} between the two sampling regimes, and the joint
CLT for the fixed-basepoint logarithm averages under local
sampling, which determines the limiting mean shifts in
Theorem~\ref{thm:moving-clts}.

\corrected{The mutual contiguity was proved above. For the local score limits, put $W(x):=\log_\mu x$. Then $P_{r,N}W=a_r/\sqrt N$ and $\operatorname{Cov}_{P_{r,N}}(W)\to\Sigma_r$. Since $dP_{r,N}/dP_r\le3/2$ for large $N$ and $P_r\|W\|^2<\infty$, the centered triangular arrays satisfy Lindeberg's condition. The multivariate triangular-array CLT, $n_r/N\to\lambda_r$, and independence across groups therefore give}
{\color{black}
\begin{equation}\label{eq:local-score-clt}
 (\sqrt N\,U_{r,N})_{r=1}^k
 \xrightarrow{\mathcal D}
 N_{(T_\mu M)^k}\!\left((a_r)_{r=1}^k,D_\lambda\right)
 \qquad\text{under }\Prb_{s,N}.
\end{equation}
}

\subsubsection{Common Hessian expansion and intrinsic cancellation: proof of Theorem~\ref{thm:moving-clts}\textup{(i)}}\label{subsec:Common Hessian expansion and intrinsic cancellation}
\begin{proof}[Proof of Theorem~\ref{thm:moving-clts}\textup{(i)}]
\revised{We first work under the baseline null model.}
Choose $\rho>0$ small enough that $\overline{B_g(\mu,3\rho)}\subset B$ and, for $U:=B_g(\mu,\rho)$, the two-point logarithm is smooth on a neighborhood of $\overline U\times\overline U$. Such a choice is available because the complement of the cut relation is open and contains $(\mu,\mu)$; no geodesic convexity is needed. For $x$ away from $\Cut(\mu)$, write $K_\mu(x):=\Hess_{p=\mu}^{\sharp}d_g(p,x)^2$, the $(1,1)$-Hessian with respect to the first argument. Since $\log_p x=-\tfrac12\operatorname{grad}_p d_g(p,x)^2$ off the cut relation, Taylor expansion along $a=\exp_\mu u$, with parallel transport back to $\mu$, gives
\begin{equation}\label{eq:common-hessian-expansion}
\tau_{a,\mu}\log_a x=\log_\mu x-\tfrac12K_\mu(x)\log_\mu a+R(a,x),\qquad \|R(a,x)\|\le C_0d_g(a,\mu)^2.
\end{equation}
The bound is uniform for $a,x\in U$: the transported logarithm is smooth on a neighborhood of $\overline U\times\overline U$, so its second derivatives in $u$ are bounded. Under Assumption~\ref{as:compact-support}, the same expansion holds uniformly for $x\in C$ and $a$ in a sufficiently small fixed normal ball about $\mu$. Indeed, the complement of the cut relation is open, and compactness of $C$ supplies a product neighborhood of $\{\mu\}\times C$ on which these derivatives are bounded.

For the intrinsic statistic, substitute $x=\widehat\mu_{r,n_r}$. In normal coordinates $p=\exp_\mu v$, $d_g(p,\mu)^2=\|v\|^2$, whose first derivative at $0$ is zero and second derivative is $2I_d$; hence $K_\mu(\mu)=2I_d$. Continuity and groupwise consistency give $K_\mu(\widehat\mu_{r,n_r})\inPrb2I_d$. Consequently, if $\sqrt N\xi_N=O_p(1)$, \eqref{eq:common-hessian-expansion} yields, simultaneously for all groups,
\begin{equation}\label{eq:intrinsic-hessian-linearization}
\sqrt N\,\tau_{\widehat\mu_N,\mu}\log_{\widehat\mu_N}\widehat\mu_{r,n_r}
=\sqrt N\log_\mu\widehat\mu_{r,n_r}-\sqrt N\xi_N+o_p(1).
\end{equation}
The quadratic remainder is $O_p(N^{-1/2})$ after multiplication by $\sqrt N$, and the Hessian error times $\sqrt N\xi_N$ is $o_p(1)$. Thus the same first-order correction cancels between groups $r$ and $k$. Proposition~\ref{prop:pooled-clt} supplies this root-$N$ tightness when its hypotheses hold, but it is not needed for the stronger intrinsic assertion.

To prove that assertion using only pooled-mean consistency, put $G(a,x,y):=\tau_{a,\mu}(\log_a x-\log_a y)$. Appendix~\ref{app:intrinsic-cancellation} proves the uniform bound \eqref{eq:intrinsic-cancellation} for $a,x,y\in U$. By Lemma~\ref{lem:pooled-consistency}, all sample means lie in $U$ with probability tending to one. The groupwise CLTs give $\sqrt N\,d_g(\widehat\mu_{r,n_r},\widehat\mu_{k,n_k})=O_p(1)$, whereas pooled consistency gives $d_g(\widehat\mu_N,\mu)=o_p(1)$. Therefore, for every $r<k$,
\begin{equation}\label{eqn:diff-PT-log-intrinsic}
   \sqrt N\left\|
 \tau_{\widehat\mu_N,\mu}\Delta_{r,N}^{(i)}
 -(\log_\mu\widehat\mu_{r,n_r}-\log_\mu\widehat\mu_{k,n_k})
 \right\|
 \le C d_g(\widehat\mu_N,\mu)\,
 \sqrt N\,d_g(\widehat\mu_{r,n_r},\widehat\mu_{k,n_k})=o_p(1).  
\end{equation}

Since $k$ is fixed, $\sqrt N(\mathcal T_N\Delta_N^{(i)}-\Delta_{\mu,N}^{(i)})=o_p(1)$ in the product tangent space. Independence and \eqref{eq:groupwise-clt} give independent Gaussian limits $Y_r$ with covariance $Q_r/\lambda_r$ for $(\sqrt N\log_\mu\widehat\mu_{r,n_r})_{r=1}^k$. Hence the contrast limit is $(Y_r-Y_k)_{r<k}$, with blocks $\delta_{rs}Q_r/\lambda_r+Q_k/\lambda_k$, including the off-diagonal blocks. Slutsky's theorem proves \eqref{eq:intrinsic-moving-clt}.

\revised{\textit{Extension under local sampling.} The groupwise M-estimator argument used for \eqref{eq:groupwise-clt} also gives $\sqrt N\log_\mu\widehat\mu_{r,n_r}=2H_r^{-1}\sqrt N\,U_{r,N}+o_p(1)$ under $\Prb_{0,N}$, after multiplying the groupwise linearization by $\sqrt{N/n_r}$. Contiguity transfers this remainder and pooled-mean consistency to $\Prb_{s,N}$. The shifted score limits above therefore give independent groupwise Gaussian limits with means $2H_r^{-1}a_r$ and covariances $Q_r/\lambda_r$. The deterministic bound \eqref{eq:intrinsic-cancellation} in Appendix~\ref{app:intrinsic-cancellation}, or equivalently transfer of the remainder in \eqref{eqn:diff-PT-log-intrinsic}, gives the same $o_p(1)$ difference between the scaled moving and fixed contrasts under $\Prb_{s,N}$. Slutsky's theorem proves \eqref{eq:intrinsic-moving-clt} with mean $h^{(i)}$ and the unchanged covariance blocks $\delta_{rs}Q_r/\lambda_r+Q_k/\lambda_k$. This argument uses only pooled consistency, not a pooled-mean CLT or either alternative of Lemma~\ref{lem:empirical-moving-log}.}
\end{proof}
\subsubsection{Empirical linearization and arithmetic correction: proof of Theorem~\ref{thm:moving-clts}\textup{(ii)}}\label{sec:moving-log-linearization}
For the arithmetic statistic, put 
\begin{enumerate}
    \item the parallel transport of the logarithm vector field $h_a(x):=\tau_{a,\mu}\log_a x\in T_\mu M$, so $h_\mu=\log_\mu$.
    \item the expected logarithm vector field $G_r(a):=\int_M\log_a x\,dP_r(x)\in T_aM$.
    \item For each group $r$, let the empirical measure be $P_{r,n_r}:=n_r^{-1}\sum_{j=1}^{n_r}\delta_{X_j^{(r)}}$. Following \citet{vanDerVaart1998}, use measures as \textit{integral operators}: $P_rf:=\int_M f(x)\,dP_r(x)$ and $P_{r,n_r}f:=\int_M f(x)\,dP_{r,n_r}(x)=n_r^{-1}\sum_{j=1}^{n_r}f(X_j^{(r)})$.
\end{enumerate}

\begin{lem}[\textbf{Empirical moving-logarithm linearization}]\label{lem:empirical-moving-log}
Assume $H_0$ \revised{or the local sampling regime of Theorem~\ref{thm:moving-clts}} and $n_r/N\to\lambda_r\in(0,1)$. Suppose either \textup{(A)} Assumption~\ref{as:compact-support} holds, or \textup{(B)} $M$ is compact, connected and without boundary, and each $P_r$ satisfies assumptions \textup{(iii)--(v)} of \citet[Theorem~1]{HotzLeWood2024}, with the common population Fréchet mean $\mu$ and the probability measure in \citet{HotzLeWood2024}
replaced by $P_r$. In case \textup{(B)}, assume additionally that, for every $r$ and $\varepsilon>0$, the following stochastic equicontinuity condition at $\mu$:
\begin{equation}\label{eq:moving-log-equicontinuity}
\lim_{\delta\downarrow0}\limsup_{N\to\infty}\Prb\!\left(\sup_{a\in B_\delta(\mu)}\left\|\sqrt{n_r}(P_{r,n_r}-P_r)(h_a-h_\mu)\right\|>\varepsilon\right)=0.
\end{equation}
Here $B_\delta(\mu)$ lies in a fixed normal neighborhood, the suprema are assumed measurable, and the probabilities are ordinary probabilities. Condition~\eqref{eq:moving-log-equicontinuity} is local asymptotic stochastic equicontinuity at $\mu$ of the groupwise empirical process $\mathbb G_{r,n_r}(f):=\sqrt{n_r}(P_{r,n_r}-P_r)f$, indexed by $h_a$.

\revised{Under local sampling, the $P_r$ remain the baseline laws, the observations have laws $P_{r,N}$, and \eqref{eq:moving-log-equicontinuity} is assumed under $\Prb_{0,N}$. All convergence and stochastic orders for $a_N$ and the conclusion below are taken under the sampling regime under consideration.}

Let $a_N$ be any measurable, possibly data-dependent sequence with $a_N\to\mu$ in probability and $\sqrt N\log_\mu a_N=O_p(1)$, where the logarithm is understood on the event that $a_N$ lies in a fixed normal neighborhood of $\mu$. Then, simultaneously for $r=1,\ldots,k$,
\begin{equation}\label{eq:empirical-moving-log}
\tau_{a_N,\mu}\left(\frac1{n_r}\sum_{j=1}^{n_r}\log_{a_N}(X_j^{(r)})\right)=U_{r,N}-\frac12H_r\log_\mu a_N+o_p(N^{-1/2}).
\end{equation}
Here $H_r=\Hess_\mu^\sharp F_r$ is the full population Fr\'echet Hessian, including the cut-locus contribution when present.
\end{lem}

For the proof, please see Appendix \ref{app:empirical-moving-log}.

\begin{proof}[Proof of Theorem~\ref{thm:moving-clts}\textup{(ii)}]
\revised{We first work under the baseline null model.}
Put $\xi_N:=\log_\mu\widehat\mu_N$. Lemma~\ref{lem:pooled-consistency} and Proposition~\ref{prop:pooled-clt} give $\widehat\mu_N\to\mu$ in probability and $\sqrt N\xi_N=O_p(1)$. Apply Lemma~\ref{lem:empirical-moving-log}, whose proof is in Appendix~\ref{app:empirical-moving-log}, with $a_N=\widehat\mu_N$ to obtain
\begin{equation}\label{eq:arithmetic-linearization}
\sqrt N\,\tau_{\widehat\mu_N,\mu}V^{(a)}_{r,N}
=
\sqrt N\,U_{r,N}
-\frac12H_r\sqrt N\,\xi_N+o_p(1).
\end{equation}

Under alternative \textup{(A)}, this is the average of \eqref{eq:common-hessian-expansion}: averaging its uniformly bounded remainder gives an $O_p(N^{-1/2})$ error at the $\sqrt N$ scale, and $P_{r,n_r}K_\mu\inPrb H_r$ by the law of large numbers and differentiation under the integral. Under alternative \textup{(B)}, the population expansion \eqref{eq:population-moving-log} and the explicit empirical equicontinuity condition \eqref{eq:moving-log-equicontinuity} yield the same linearization with the full $H_r$, including its cut-locus contribution; pointwise Hessian averaging is not used.

Subtracting group $k$ from group $r$ gives $\sqrt N(U_{r,N}-U_{k,N})-\tfrac12(H_r-H_k)\sqrt N\xi_N+o_p(1)$. Unlike the intrinsic correction, this generally survives. Substituting the joint linear representation \eqref{eq:pooled-linearization} gives
\begin{equation}\label{eq:arithmetic-contrast}
\sqrt N\,\tau_{\widehat\mu_N,\mu}
(V^{(a)}_{r,N}-V^{(a)}_{k,N})
=
\sum_{s=1}^k B_{rs}\sqrt N\,U_{s,N}
+o_p(1),
\quad 1\leq r<k.
\end{equation}

Replacing $\lambda_{s,N}$ by $\lambda_s$ costs $o_p(1)$ because $k$ is fixed and every $\sqrt N U_{s,N}=O_p(1)$. The independent groupwise multivariate CLTs give $(\sqrt N U_{s,N})_{s=1}^k\xrightarrow{\mathcal D}N_{(T_\mu M)^k}(0,D_\lambda)$; hence \eqref{eq:arithmetic-contrast} has covariance $BD_\lambda B^*$. Both the fixed-basepoint term and $\sqrt N\xi_N$ depend on these same averages and are generally correlated; no independence between them is assumed. This proves \eqref{eq:arithmetic-moving-clt}. \corrected{Applying \citet[Theorem~1.4]{LytchakHuckemann2025} to the pooled empirical measure, whose Fr\'echet mean is $\widehat\mu_N$, gives $X_j^{(r)}\notin\Cut(\widehat\mu_N)$ for every $r,j$. Hence the ordinary and zero-extended statistics coincide and have the same limit.}

\revised{\textit{Extension under local sampling.} Proposition~\ref{prop:pooled-clt} proves that the difference between the two sides of \eqref{eq:pooled-linearization}, excluding its $o_p(1)$ term, tends to zero in probability under $\Prb_{0,N}$. This difference is a measurable function of the same observations and estimators, so contiguity makes it $o_p(1)$ under $\Prb_{s,N}$ as well. In particular, the shifted score limits imply $\sqrt N\xi_N=O_p(1)$ under $\Prb_{s,N}$; pooled consistency transfers in the same way. Lemma~\ref{lem:empirical-moving-log} therefore applies with $a_N=\widehat\mu_N$, and \eqref{eq:arithmetic-linearization}--\eqref{eq:arithmetic-contrast} retain their baseline Hessians and remainder orders. Replacing $\lambda_{s,N}$ by $\lambda_s$ still costs $o_p(1)$, because the shifted score limits give $\sqrt N U_{s,N}=O_p(1)$. The continuous mapping theorem and Slutsky's theorem applied to \eqref{eq:arithmetic-contrast} give mean $B(a_1,\ldots,a_k)=h^{(a)}$ and covariance $BD_\lambda B^*=\Gamma^{(a)}$. The pooled correction and fixed-basepoint averages still depend on the same scores; no independence between them is used. \corrected{The same empirical-measure argument shows that the ordinary and zero-extended statistics also coincide under local sampling.}}
\end{proof}

\begin{Remark}[\textbf{Smeary intrinsic limits are also not altered by moving basepoint}]\label{rem:smeary-moving}
The cancellation underlying Theorem~\ref{thm:moving-clts}\textup{(i)} does not rely
on the $N^{-1/2}$ rate.  Indeed following the smeariness introduced in \citet{EltznerHuckemann2019}, drop the positive-definiteness assumption on
the $H_r$ and suppose instead that, for some $\kappa>0$ and
$\alpha_\kappa:=1/\{2(\kappa+1)\}$,
$n_r^{\alpha_\kappa}\log_\mu\widehat\mu_{r,n_r}\inD Z_r$,
$r=1,\ldots,k$, where the $Z_r$ need not be Gaussian. With
$\Delta_{\mu,N}^{(i)}$ as defined above, we have
\begin{equation}
 N^{\alpha_\kappa}
 \bigl(\mathcal T_N\Delta_N^{(i)}-\Delta_{\mu,N}^{(i)}\bigr)\inPrb0.       \label{eq:smeary-equivalence}
\end{equation}
Moreover, since the groups are independent and $n_r/N\to\lambda_r$,
\[
 N^{\alpha_\kappa}\Delta_{\mu,N}^{(i)}
 \inD
 \bigl(\lambda_1^{-\alpha_\kappa}Z_1-\lambda_k^{-\alpha_\kappa}Z_k,
 \ldots,
 \lambda_{k-1}^{-\alpha_\kappa}Z_{k-1}
       -\lambda_k^{-\alpha_\kappa}Z_k\bigr),
\]
and the same limit therefore holds for
$N^{\alpha_\kappa}\mathcal T_N\Delta_N^{(i)}$.

Indeed, convergence in distribution implies
$n_r^{\alpha_\kappa}
d(\widehat\mu_{r,n_r},\mu)
=n_r^{\alpha_\kappa}\|\log_\mu\widehat\mu_{r,n_r}\|=O_p(1)$.
Since $n_r/N\to\lambda_r>0$, this gives
$N^{\alpha_\kappa}d(\widehat\mu_{r,n_r},\mu)=O_p(1)$ and hence, by the
triangle inequality,
$N^{\alpha_\kappa}
d(\widehat\mu_{r,n_r},\widehat\mu_{k,n_k})=O_p(1)$.
The proof of Lemma~\ref{lem:pooled-consistency} does not use positive
definiteness of the $H_r$, so
$d(\widehat\mu_N,\mu)=o_p(1)$ still holds.  Applying \eqref{eq:intrinsic-cancellation}, for each
$r<k$,
\[
 N^{\alpha_\kappa}
 \bigl\|G(\widehat\mu_N,\widehat\mu_{r,n_r},\widehat\mu_{k,n_k})
       -G(\mu,\widehat\mu_{r,n_r},\widehat\mu_{k,n_k})\bigr\|
 \, \le
 C\,d(\widehat\mu_N,\mu)\,
 N^{\alpha_\kappa}
 d(\widehat\mu_{r,n_r},\widehat\mu_{k,n_k})
 =o_p(1).
\]
These are precisely the $k-1$ components of
$N^{\alpha_\kappa}
(\mathcal T_N\Delta_N^{(i)}-\Delta_{\mu,N}^{(i)})$; since $k$ is fixed,
\eqref{eq:smeary-equivalence} follows.
\end{Remark}
\begin{Remark}[\textbf{Euclidean case}]
When $M=\mathbb R^d$ is equipped with the standard Euclidean
metric, parallel transport is the identity and $H_r=2I_d$
for every group $r$. Consequently, $H_r-H_k=0$, and the
additional Hessian-dependent term in \eqref{eq:arithmetic-blocks}
vanishes. In fact,
$V^{(a)}_{r,N}-V^{(a)}_{k,N}=U_{r,N}-U_{k,N}$
holds exactly, since the pooled sample mean cancels from
the difference. Thus, the arithmetic moving-basepoint CLT
reduces to the classical Euclidean contrast CLT, without
any additional contribution from the random basepoint.
Although using the pooled sample mean as the basepoint
generally changes the limiting distribution of each
individual group average, it does not affect the
limiting distribution of their contrasts.
\end{Remark}

\begin{Remark}[\textbf{Special case of two groups ($k=2$)}]\label{rem:intrinsic-two-groups}\label{rem:two-groups}
Write $\lambda_1=\lambda$ and $\lambda_2=1-\lambda$. The intrinsic conclusion of Theorem~\ref{thm:moving-clts}\textup{(i)} has covariance
\begin{equation}
\Gamma^{(i)}
=\frac4\lambda H_1^{-1}\Sigma_1H_1^{-1}
+\frac4{1-\lambda}H_2^{-1}\Sigma_2H_2^{-1}.
\label{eq:intrinsic-two-group-covariance}
\end{equation}
This is exactly the fixed-basepoint intrinsic covariance. For arithmetic contrasts, $H_\lambda=\lambda H_1+(1-\lambda)H_2$ gives $B_{11}=H_2H_\lambda^{-1}$ and $B_{12}=-H_1H_\lambda^{-1}$, so Theorem~\ref{thm:moving-clts}\textup{(ii)} has covariance
\begin{equation}\label{eq:arithmetic-two-group-covariance}
\Gamma^{(a)}
=\frac1\lambda
 H_2H_\lambda^{-1}\Sigma_1H_\lambda^{-1}H_2
+\frac1{1-\lambda}
 H_1H_\lambda^{-1}\Sigma_2H_\lambda^{-1}H_1.
\end{equation}
The order of the factors is retained; the Hessians need not commute. If $H_1=H_2$, this reduces to the fixed-basepoint covariance $\Sigma_1/\lambda+\Sigma_2/(1-\lambda)$.
\end{Remark}

\section{Plug-in covariance estimation and the Wald consequence}\label{sec:wald}

We develop Wald-type tests based on the intrinsic and arithmetic moving-basepoint CLTs, Theorem~\ref{thm:moving-clts}\textup{(i)--(ii)}, respectively. Throughout this section, assume Assumption~\ref{as:compact-support} and the groupwise conditions stated for \eqref{eq:groupwise-clt}, including
positive definiteness of every $H_r$.

\subsection[Wald-type tests based on the intrinsic and arithmetic moving-basepoint CLTs, Theorem \ref{thm:moving-clts}\textup{(i)--(ii)}]{Wald-type tests based on the intrinsic and arithmetic\\ moving-basepoint CLTs, Theorem \ref{thm:moving-clts}\textup{(i)--(ii)}}
\label{subsec:wald-tests}

\subsubsection{Plug-in covariance estimation and studentization}
\label{subsec:studentization}

For $r=1,\ldots,k$, set
\[
 Z_{rj,N}:=\log_{\widehat\mu_N}X^{(r)}_j,
 \qquad
 \overline Z_{r,N}:=\frac1{n_r}\sum_{j=1}^{n_r}Z_{rj,N},
\]
and define on $T_{\widehat\mu_N}M$ the covariance operator on the moving tangent space $T_{\widehat\mu_N}M$ (N.B: it is not intrinsically a matrix unless we choose a basis)
\[
 \widehat\Sigma_{r,N}v:=\frac1{n_r-1}\sum_{j=1}^{n_r}
 \langle Z_{rj,N}-\overline Z_{r,N},v\rangle
 (Z_{rj,N}-\overline Z_{r,N}).
\]

\subsubsection{Plug-in Hessian estimation and studentization, and why Assumption \ref{as:compact-support} is used here}

Let $\widehat H_{r,N}:T_{\widehat\mu_N}M\to
T_{\widehat\mu_N}M$ be the self-adjoint endomorphism
determined by
$\langle\widehat H_{r,N}u,v\rangle_{\widehat\mu_N}
=\Hess_{\widehat\mu_N}F_{r,n_r}(u,v)$.
Under Assumption~\ref{as:compact-support}, the
squared-distance functions are smooth at
$\widehat\mu_N$ with probability tending to one,
and therefore
\begin{equation}\label{eq:empirical-hessian-operator}
\widehat H_{r,N}
=\frac1{n_r}\sum_{j=1}^{n_r}
\Hess_{\widehat\mu_N}^{\sharp}
d_g(\cdot,X_j^{(r)})^2.
\end{equation}

We retain Assumption~\ref{as:compact-support} for the plug-in estimation because, in its absence, the full population Hessian $H_r$ may differ from the expected pointwise squared-distance Hessian
by the cut-locus correction $-2J_{P_r}(\mu)$;
see \citet[Theorem~1, equations~(14)--(15)] {HotzLeWood2024}.
The correction $J_{P_r}(\mu)$ is given explicitly by a surface integral over $\mathcal H_\mu$ in \citet[equation~(20), p.~1228]{HotzLeWood2024}. By \citet[Lemma~1, p.~1226]{HotzLeWood2024},
$\mathcal H_\mu$ is a countable union of disjoint codimension-one submanifolds, each of whose points is joined to $\mu$ by exactly two minimizing geodesics. Consistently estimating this correction
requires identifying the relevant cut-locus geometry and estimating the probability density and geometric factors appearing in the surface
integral. This is generally nontrivial, and the empirical Hessian in
\eqref{eq:empirical-hessian-operator} need not
consistently estimate the full population Hessian
when the cut-locus correction is nonzero.

Whenever the required inverses exist, put
$\widehat H_{\lambda,N}:=\sum_r\lambda_{r,N}
\widehat H_{r,N}$ and define
$\widehat\Gamma_N^{(i)}$ from
\eqref{eq:intrinsic-covariance} by replacing
$(\lambda_r,H_r,\Sigma_r)$ with
$(\lambda_{r,N},\widehat H_{r,N},
\widehat\Sigma_{r,N})$.
Likewise define $\widehat B_N$, $\widehat D_N$
and
$\widehat\Gamma_N^{(a)}
:=\widehat B_N\widehat D_N\widehat B_N^*$
from \eqref{eq:arithmetic-blocks}--
\eqref{eq:arithmetic-covariance}
by the same replacements.
Set these estimators equal to zero when the
required inverses do not exist.

\begin{lem}[\textbf{Consistency of the plug-in covariance operators}]\label{lem:cov-consistency}
Under the standing assumptions of this section, for each $\star\in\{i,a\}$,
\begin{equation}
\mathcal T_N\widehat\Gamma_N^{(\star)}\mathcal T_N^*\xrightarrow{\mathrm{a.s.}}\Gamma^{(\star)}. \label{eq:cov-consistency}
\end{equation}
\revised{Under the local sampling regime of Theorem~\ref{thm:moving-clts}, the same limit holds in probability. Also, writing $\tau_N:=\tau_{\widehat\mu_N,\mu}$, the groupwise consistency $\tau_N\widehat\Sigma_{r,N}\tau_N^*\inPrb\Sigma_r$ transfers under the weaker baseline hypotheses of Lemma~\ref{lem:pooled-consistency} and Assumption~\ref{as:compact-support}; it does not require positive definiteness of each $H_r$.}
\end{lem}

\begin{proof}
Write $\tau_N:=\tau_{\widehat\mu_N,\mu}$, so that $\mathcal T_N$ acts componentwise as $\tau_N$, and let $C:=\bigcup_{r=1}^k\operatorname{supp}(P_r)$. By Assumption \ref{as:compact-support}, $C$ is compact and separated from $\operatorname{Cut}_M(\mu)$. Hence, on a sufficiently small neighborhood of $\mu$, the maps $(a,x)\mapsto\tau_{a,\mu}\log_a x$ and $(a,x)\mapsto\tau_{a,\mu} \circ \operatorname{Hess}^{\sharp}_{p=a}d_g(p,x)^2 \circ \tau_{a,\mu}^*$ are jointly smooth, with uniformly bounded derivatives in $a$ for $x\in C$.

By Lemma~\ref{lem:pooled-consistency}, $\widehat\mu_N\to\mu$ almost surely. Consequently, for some deterministic $K<\infty$, almost surely for all sufficiently large $N$,
\[
\varepsilon_N:=\sup_{x\in C}\|\tau_N\log_{\widehat\mu_N}x-\log_\mu x\|\le Kd_g(\widehat\mu_N,\mu)\longrightarrow0.
\]
For each group $r$, put $Y_{rj,N}:=\tau_N\log_{\widehat\mu_N}X_j^{(r)}$ and $W_{rj}:=\log_\mu X_j^{(r)}$, and let $\widehat\Sigma_{r,N}^{0}$ be the ordinary sample covariance of the $W_{rj}$ in $T_\mu M$. Since $\tau_N$ is a linear isometry, with $\tau_N^*=\tau_N^{-1}$, the operator $\tau_N\widehat\Sigma_{r,N}\tau_N^*$ is exactly the sample covariance of the $Y_{rj,N}$. Writing $L:=\sup_{x\in C}\|\log_\mu x\|<\infty$ and using $\|u\otimes u-v\otimes v\|_{\mathrm{op}}\le(\|u\|+\|v\|)\|u-v\|$, we obtain
\[
\|\tau_N\widehat\Sigma_{r,N}\tau_N^*-\widehat\Sigma_{r,N}^{0}\|_{\mathrm{op}}\le\frac{2n_r}{n_r-1}(2L+\varepsilon_N)\varepsilon_N\longrightarrow0\quad\text{a.s.}
\]
The ordinary strong law applied to $W_{rj}$ and $W_{rj}\otimes W_{rj}$ gives $\widehat\Sigma_{r,N}^{0}\to\Sigma_r$ almost surely. Thus $\tau_N\widehat\Sigma_{r,N}\tau_N^*\to\Sigma_r$ almost surely, without requiring independence of the $Y_{rj,N}$.

Similarly, writing $K(a,x):=\operatorname{Hess}^{\sharp}_{p=a}d_g(p,x)^2$, uniform smoothness and $\widehat\mu_N\to\mu$ almost surely give $\sup_{x\in C}\|\tau_NK(\widehat\mu_N,x)\tau_N^*-K(\mu,x)\|_{\mathrm{op}}\to0$ almost surely. Assumption \ref{as:compact-support} permits differentiation under the population integral, so $H_r=\int_M K(\mu,x)\,dP_r(x)$. The ordinary strong law therefore yields $\tau_N\widehat H_{r,N}\tau_N^*\to H_r$ almost surely, simultaneously for all $r$.

Since every $H_r$ and $H_\lambda$ is positive definite and $\lambda_{r,N}\to\lambda_r$, continuity of inversion, addition and multiplication of operators in the fixed space $T_\mu M$ proves \eqref{eq:cov-consistency}. The asserted eventual invertibility follows from the same convergence.

\revised{For local sampling, each convergence-in-probability consequence of the preceding baseline almost-sure statements transfers by contiguity. This proves the asserted local consistency of $\widehat\Gamma_N^{(\star)}$ and invertibility with probability tending to one. The groupwise covariance argument above uses only compact support away from the cut locus, baseline finite second moments and pooled consistency, before any Hessian inversion; hence its stated weaker hypotheses suffice. For $k=2$, \eqref{eq:SPN} and $n/N\to\lambda$ consequently give $\tau_NS_{p,N}\tau_N^*\inPrb\Sigma_\lambda$ and $\tau_N\widehat\Gamma_{\mathrm{pool},N}\tau_N^*\inPrb\Gamma_{\mathrm{pool}}=\Sigma_\lambda/[\lambda(1-\lambda)]$ under local sampling. These are probability limits, with no assertion of almost-sure convergence under the triangular sampling regime.}
\end{proof}

\begin{Cor}[\textbf{Asymptotic chi-square laws}]\label{cor:wald}
Assume that $\Gamma^{(\star)}$ is positive definite for $\star\in\{i,a\}$. On the event that $\widehat\Gamma_N^{(\star)}$ is invertible, define
\[
T_N^{(\star)}:=N\left\langle\Delta_N^{(\star)},(\widehat\Gamma_N^{(\star)})^{-1}\Delta_N^{(\star)}\right\rangle_{(T_{\widehat\mu_N}M)^{k-1}},
\]
and set $T_N^{(\star)}:=0$ otherwise. Then, under $H_0$, $T_N^{(\star)}\inD\chi^2_{(k-1)d}$.
\end{Cor}

\begin{proof}
Put $S_N^{(\star)}:=\sqrt N\,\mathcal T_N\Delta_N^{(\star)}$. By Theorem~\ref{thm:moving-clts}\textup{(i)} or \textup{(ii)},
\[
S_N^{(\star)}\inD N_{(T_\mu M)^{k-1}}(0,\Gamma^{(\star)}).
\]
By Lemma~\ref{lem:cov-consistency}, $\mathcal T_N\widehat\Gamma_N^{(\star)}\mathcal T_N^*\to\Gamma^{(\star)}$ almost surely and, since $\Gamma^{(\star)}$ is positive definite, $\widehat\Gamma_N^{(\star)}$ is eventually invertible almost surely. Because $\mathcal T_N$ is a linear isometry, $\mathcal T_N^*=\mathcal T_N^{-1}$ (so $\mathcal T_N^*\mathcal T_N=I$), and the \textit{quadratic form is invariant under parallel transport}:
\begin{multline}\label{eq:wald-transport-invariance}
\left\langle\Delta_N^{(\star)},(\widehat\Gamma_N^{(\star)})^{-1}\Delta_N^{(\star)}\right\rangle_{(T_{\widehat\mu_N}M)^{k-1}}
=\left\langle\mathcal T_N\Delta_N^{(\star)},\mathcal T_N(\widehat\Gamma_N^{(\star)})^{-1}\mathcal T_N^*\mathcal T_N\Delta_N^{(\star)}\right\rangle_{(T_\mu M)^{k-1}}\\
=\left\langle\mathcal T_N\Delta_N^{(\star)},(\mathcal T_N\widehat\Gamma_N^{(\star)}\mathcal T_N^*)^{-1}\mathcal T_N\Delta_N^{(\star)}\right\rangle_{(T_\mu M)^{k-1}}.
\end{multline}
Here the second equality uses $(\mathcal T_N\widehat\Gamma_N^{(\star)}\mathcal T_N^*)^{-1}=\mathcal T_N(\widehat\Gamma_N^{(\star)})^{-1}\mathcal T_N^*$.

Continuity of inversion and Slutsky's theorem, applied in the fixed tangent space $(T_\mu M)^{k-1}$, now give $T_N^{(\star)}\inD\langle Z,(\Gamma^{(\star)})^{-1}Z\rangle\sim\chi^2_{(k-1)d}$, where $Z\sim N_{(T_\mu M)^{k-1}}(0,\Gamma^{(\star)})$.

\end{proof}

\begin{Test}[\textbf{Multisample Wald-type tests for equality of Fr\'echet means}]\label{test:multisample-wald}
For $\star\in\{i,a\}$ and $0<\alpha<1$, reject the common-Fr\'echet-mean null hypothesis $H_0$ whenever
\[
T_N^{(\star)}>q_{(k-1)d,1-\alpha},
\]
where $T_N^{(\star)}$ is defined in Corollary~\ref{cor:wald} and $q_{\nu,1-\alpha}$ denotes the $(1-\alpha)$-quantile of $\chi_\nu^2$. Under the assumptions of Corollary~\ref{cor:wald}, each test has asymptotic Type~I error $\alpha$; for $k=2$, the reference distribution is $\chi_d^2$. These tests are manifold analogues of classical multisample Wald tests for equality of multivariate mean vectors; see \citet[Chapter~9, Section~9.2.4(b), pp.~445--447]{Seber1984}.
\end{Test}

\begin{Cor}[\textbf{Local asymptotic power of the corrected Wald tests}]\label{cor:wald-local-power}
Assume the standing hypotheses of Section~\ref{sec:wald} and that $\Gamma^{(\star)}$ is positive definite for $\star\in\{i,a\}$. Let the observations follow the local sampling regime of Theorem~\ref{thm:moving-clts}. Put $\nu:=(k-1)d$ and
\[
\delta_\star:=\left\langle h^{(\star)},(\Gamma^{(\star)})^{-1}h^{(\star)}\right\rangle_{(T_\mu M)^{k-1}}.
\]
Then, under the locally perturbed sampling laws, $T_N^{(\star)}\inD\chi_\nu^2(\delta_\star)$, where $\chi_\nu^2(\delta_\star)$ denotes the noncentral chi-square distribution with $\nu$ degrees of freedom and noncentrality parameter $\delta_\star$. Consequently, the asymptotic local power of the test in Test~\ref{test:multisample-wald} is
\[
\lim_{N\to\infty}\Prb_{s,N}\!\left(T_N^{(\star)}>q_{\nu,1-\alpha}\right)=1-F_{\chi_\nu^2(\delta_\star)}\!\left(q_{\nu,1-\alpha}\right).
\]
Whenever the perturbed sampling laws violate $H_0$, the corresponding asymptotic Type~II error is $F_{\chi_\nu^2(\delta_\star)}(q_{\nu,1-\alpha})$. In particular, $\delta_\star=0$ gives limiting rejection probability $\alpha$, but does not necessarily imply that the perturbed laws satisfy $H_0$.
\end{Cor}

\begin{proof}
Put $S_N^{(\star)}:=\sqrt N\,\mathcal T_N\Delta_N^{(\star)}$. By Theorem~\ref{thm:moving-clts}, under the local sampling regime, $S_N^{(\star)}\inD Z_\star$, where $Z_\star\sim N_{(T_\mu M)^{k-1}}(h^{(\star)},\Gamma^{(\star)})$. By Lemma~\ref{lem:cov-consistency}, $\mathcal T_N\widehat\Gamma_N^{(\star)}\mathcal T_N^*\inPrb\Gamma^{(\star)}$ under the same regime. Since $\Gamma^{(\star)}$ is positive definite, $\widehat\Gamma_N^{(\star)}$ is invertible with probability tending to one. The parallel-transport invariance in \eqref{eq:wald-transport-invariance}, continuity of inversion and Slutsky's theorem therefore give
\[
T_N^{(\star)}\inD\left\langle Z_\star,(\Gamma^{(\star)})^{-1}Z_\star\right\rangle.
\]
To identify the limiting distribution, choose an orthonormal basis of $(T_\mu M)^{k-1}$ and put $W_\star:=(\Gamma^{(\star)})^{-1/2}Z_\star$. Then
\[
W_\star\sim N_\nu\!\left((\Gamma^{(\star)})^{-1/2}h^{(\star)},I_\nu\right),\qquad \nu:=(k-1)d.
\]
Thus the coordinates of $W_\star$ are independent normal random variables with unit variances, and
\[
\left\langle Z_\star,(\Gamma^{(\star)})^{-1}Z_\star\right\rangle=\|W_\star\|^2\sim\chi_\nu^2(\delta_\star),
\]
where
\[
\delta_\star=\left\|(\Gamma^{(\star)})^{-1/2}h^{(\star)}\right\|^2=\left\langle h^{(\star)},(\Gamma^{(\star)})^{-1}h^{(\star)}\right\rangle.
\]
Since the noncentral chi-square distribution has a continuous distribution function, the asserted limiting rejection and non-rejection probabilities follow.
\end{proof}

\subsection{Comparison with an earlier pooled-sample mean Hotelling test, with suggested correction}
\label{subsec:huckemann-hotelling}

\citet[Section~6.1, Theorem~10, pp.~1244--1245]{Huckemann2012}
consider a two-sample Hotelling statistic formed from individual
coordinates $\phi_{\widehat\mu_N}(X_j^{(r)})$ at a pooled sample
Fr\'echet mean. Their admissible charts include intrinsic tangent-space
coordinates $\phi_p(x)=\log_p x$. Under the common-mean null, the theorem
asserts that the resulting statistic is ``asymptotically Hotelling
$T^2$-distributed'' if the sample sizes are asymptotically equal or
the fixed-pole covariance matrices are equal. Since the finite-sample
Hotelling distribution depends on sample size and dimension, this
phrase does not specify a single limiting law. The natural classical
interpretation for a fixed tangent-space dimension $d$ is convergence
to $\chi_d^2$: indeed, under the usual Gaussian equal-covariance model,
the two-sample statistic $T^2$ satisfies
\[
 \frac{N-d-1}{d(N-2)}T^2\sim F_{d,N-d-1},\qquad N=n+m,
\]
and hence $T^2\inD\chi_d^2$ as $N\to\infty$.
For the fixed-basepoint logarithmic observations, the same chi-squared limit
follows from the ordinary multivariate CLT and covariance consistency
under either of the two conditions above, without normality.

We specialize Theorem~10 of \citet{Huckemann2012} to the intrinsic Fréchet mean and Riemannian logarithmic coordinates, that is, in their notation, \(\varphi_p(x):=\log_p x\). Write \(n=n_1\), \(m=n_2, N=n+m\), and \(n/N\to\lambda\in(0,1)\).
Let $S_{p,N}$ be the ordinary within-group pooled sample covariance matrix of the coordinates $Z_{rj,N}:=\log_{\widehat\mu_N}X_j^{(r)}$, represented in any orthonormal basis of $T_{\widehat\mu_N}M$, explicitly given by

\begin{equation}\label{eq:SPN}
\begin{aligned}
S_{p,N}&=\frac{(n-1)\widehat{\Sigma}_{1,N}+(m-1)\widehat{\Sigma}_{2,N}}{n+m-2},\\
\widehat{\Sigma}_{r,N}&=\frac{1}{n_r-1}\sum_{j=1}^{n_r}
\left(Z_{rj,N}- \overline Z_{r,N}\right)
\left(Z_{rj,N}- \overline Z_{r,N}\right)^{\top},\quad r=1,2.
\end{aligned}
\end{equation}

where $\overline Z_{r,N}:=n_r^{-1}\sum_{j=1}^{n_r}Z_{rj,N}$, $n_1=n$, $n_2=m$, and $N=n+m$. Define the classical two-sample statistic
\[
 T_{\mathrm{pool},N}^2
 :=\frac{nm}{N}\left\langle
 V_{1,N}^{(a)}-V_{2,N}^{(a)},
 S_{p,N}^{-1}(V_{1,N}^{(a)}-V_{2,N}^{(a)})
 \right\rangle,
\]
whenever $S_{p,N}$ is invertible. Under Assumption~\ref{as:compact-support},
The groupwise covariance consistency established in the proof of
Lemma~\ref{lem:cov-consistency} therefore yields
\(
\tau_NS_{p,N}\tau_N^*
\inPrb\lambda\Sigma_1+(1-\lambda)\Sigma_2
=:\Sigma_\lambda.
\)

If $\lambda=1/2$ or $\Sigma_1=\Sigma_2$ as in \citet{Huckemann2012}, Theorem 10, then $\Sigma_\lambda=\lambda(1-\lambda)\Gamma_{\mathrm{fix}}$, where $\Gamma_{\mathrm{fix}}:=\Sigma_1/\lambda+\Sigma_2/(1-\lambda)$. Consequently, since $nm/N^2\to\lambda(1-\lambda)$, the rescaled pooled sample covariance satisfies
\(
\frac{N^2}{nm}\tau_N S_{p,N}\tau_N^*\xrightarrow{\mathbb P}\Gamma_{\mathrm{fix}}.
\)
Thus, the classical pooled-covariance normalization consistently estimates the \textit{fixed-basepoint covariance} $\Gamma_{\mathrm{fix}}$, but not necessarily the \textit{moving-basepoint covariance} $\Gamma^{(a)}$.
However, Theorem~\ref{thm:moving-clts}\textup{(ii)} and
Remark~\ref{rem:two-groups} give the \emph{moving-basepoint}
covariance $\Gamma^{(a)}$, which generally differs from
$\Gamma_{\mathrm{fix}}$.
Consequently, when $\Gamma_{\mathrm{fix}}$ is positive definite, thanks to Proposition \ref{prop:pooled-hotelling-limit} or its Corollary \ref{cor:pooled-hotelling-huckemann},
\(
 T_{\mathrm{pool},N}^2 \xrightarrow{\mathcal D}
 \sum_{j=1}^{d}\rho_j\chi_{1,j}^2,
 \label{eq:uncorrected-hotelling-limit}
\)
under either of the two stated conditions, where
$\chi_{1,1}^2,\ldots,\chi_{1,d}^2$ are independent chi-squared
random variables with one degree of freedom, and
$\rho_1,\ldots,\rho_d$ are the eigenvalues of
$\Gamma_{\mathrm{fix}}^{-1/2}
\Gamma^{(a)}\Gamma_{\mathrm{fix}}^{-1/2}$.
In particular, this limiting distribution equals $\chi_d^2$
if and only if $\Gamma^{(a)}=\Gamma_{\mathrm{fix}}$.

\begin{Remark}[\textbf{Magnitude and vanishing of the covariance discrepancy}]\label{rem:covariance-equality}
For $k=2$, put the Hessian difference $D:=H_1-H_2$, $C:=DH_\lambda^{-1}$ and the tangential covariance difference $E:=\Sigma_2-\Sigma_1$, where $H_\lambda=\lambda H_1+(1-\lambda)H_2$. Expanding \eqref{eq:arithmetic-two-group-covariance} gives
\[
\Gamma^{(a)}-\Gamma_{\mathrm{fix}}=CE+EC^\top+C\Sigma_\lambda C^\top,
\]
where $\Sigma_\lambda=\lambda\Sigma_1+(1-\lambda)\Sigma_2$. Since $C^\top=H_\lambda^{-1}D$ need not equal $C$, the two linear terms cannot generally be combined. The triangle and reverse triangle inequalities yield
\[
\begin{aligned}
\|\Gamma^{(a)}-\Gamma_{\mathrm{fix}}\|_{\mathrm{op}}&\le\frac{2\|E\|_{\mathrm{op}}}{\lambda_{\min}(H_\lambda)}\|D\|_{\mathrm{op}}+\frac{\|\Sigma_\lambda\|_{\mathrm{op}}}{\lambda_{\min}(H_\lambda)^2}\|D\|_{\mathrm{op}}^2,\\
\|\Gamma^{(a)}-\Gamma_{\mathrm{fix}}\|_{\mathrm{op}}&\ge\frac{\lambda_{\min}(\Sigma_\lambda)}{\|H_\lambda\|_{\mathrm{op}}^2}\|D\|_{\mathrm{op}}^2-\frac{2\|E\|_{\mathrm{op}}}{\lambda_{\min}(H_\lambda)}\|D\|_{\mathrm{op}}.
\end{aligned}
\]
The lower bound may be negative, in which case it can be replaced by zero. In particular, if $\Sigma_1=\Sigma_2=:\Sigma$, the linear term $E$ above vanishes and $\Gamma^{(a)}-\Gamma_{\mathrm{fix}}=C\Sigma_\lambda C^\top= DH_\lambda^{-1}\Sigma H_\lambda^{-1}D\succeq0$. Moreover,
\[
\frac{\lambda_{\min}(\Sigma)}{\|H_\lambda\|_{\mathrm{op}}^2}\|D\|_{\mathrm{op}}^2
\le\|\Gamma^{(a)}-\Gamma_{\mathrm{fix}}\|_{\mathrm{op}}
=\|\Sigma^{1/2}H_\lambda^{-1}D\|_{\mathrm{op}}^2
\le\frac{\|\Sigma\|_{\mathrm{op}}}{\lambda_{\min}(H_\lambda)^2}\|D\|_{\mathrm{op}}^2.
\]
Thus, in the case of equal group covariance $\Sigma_1=\Sigma_2=:\Sigma$, the covariance discrepancy is positive semidefinite and its operator norm admits a quadratic upper bound in $\|H_1-H_2\|_{\mathrm{op}}$. If, moreover, $\Sigma$ is positive definite, the preceding lower bound is also strictly positive whenever $H_1\neq H_2$. Hence, the discrepancy admits upper and lower bounds of the same quadratic order, with constants depending on $\Sigma$ and the bounds on $H_\lambda$. Indeed, since $C\Sigma C^\top=(C\Sigma^{1/2})(C\Sigma^{1/2})^\top$, positive definiteness of $\Sigma$ implies that $C\Sigma C^\top=0$ if and only if $C=0$. As $H_\lambda$ is invertible, this yields $\Gamma^{(a)}=\Gamma_{\mathrm{fix}}$ if and only if $H_1=H_2$. If $\Sigma$ is singular, however, $C\Sigma C^\top=0$ only implies $C\Sigma^{1/2}=0$, which may hold even when $H_1\neq H_2$.

For unequal group covariances, $H_1=H_2$ remains sufficient for $\Gamma^{(a)}=\Gamma_{\mathrm{fix}}$, but need not be necessary: the terms involving $E$ may cancel the positive semidefinite term $C\Sigma_\lambda C^\top$. A complementary lower bound is available when $E$ is positive or negative definite:
\[
\|\Gamma^{(a)}-\Gamma_{\mathrm{fix}}\|_{\mathrm{op}}
\ge\frac{2\lambda_{\min}(H_\lambda)\lambda_{\min}(|E|)}{\|H_\lambda\|_{\mathrm{op}}^2}\|D\|_{\mathrm{op}}
-\frac{\|\Sigma_\lambda\|_{\mathrm{op}}}{\lambda_{\min}(H_\lambda)^2}\|D\|_{\mathrm{op}}^2,
\]
where $|E|:=(E^2)^{1/2}$. This follows by applying a Rayleigh-quotient bound to the linear term after conjugation by $H_\lambda^{-1/2}$. For fixed $H_\lambda$ and $E$, this lower bound is positive and of linear order for sufficiently small nonzero $\|D\|_{\mathrm{op}}$, whereas the preceding general lower bound can be positive when the quadratic term dominates. No commutativity assumption on $H_1$ and $H_2$ is required.

These bounds concern the discrepancy between $\Gamma^{(a)}$ and $\Gamma_{\mathrm{fix}}$. For the classical pooled Hotelling test under general unbalanced and heteroscedastic sampling, the covariance targeted by the pooled normalization is instead $\Sigma_\lambda/[\lambda(1-\lambda)]$.
\end{Remark}


\subsubsection{The inferential cost of pooled-covariance misspecification}

The distinction from the corrected Wald test is the covariance estimator, not the prefactor: defining $\widehat\Gamma_{\mathrm{pool},N}:=N^2S_{p,N}/(nm)$ gives $T_{\mathrm{pool},N}^2=N\langle\Delta_N^{(a)},\widehat\Gamma_{\mathrm{pool},N}^{-1}\Delta_N^{(a)}\rangle$. By Proposition~\ref{prop:pooled-hotelling-limit}, the rescaled pooled covariance targets $$\boxed{\Gamma_{\mathrm{pool}}:=\Sigma_\lambda/[\lambda(1-\lambda)]},$$ which equals $\Gamma_{\mathrm{fix}}$ when $\lambda=1/2$ or $\Sigma_1=\Sigma_2$. The difficulty with \citet[Theorem~10]{Huckemann2012} occurs at its final step: an $O_p(n^{-1/2})$ comparison of individual fixed- and moving-basepoint coordinates does not imply an $o_p(1)$ difference between their $\sqrt N$-scaled contrasts. Indeed, \eqref{eq:arithmetic-linearization} gives
\begin{equation}\label{eq:huckemann-missing-term}
\sqrt N\,\tau_{\widehat\mu_N,\mu}(V_{1,N}^{(a)}-V_{2,N}^{(a)})=\sqrt N(U_{1,N}-U_{2,N})-\tfrac12(H_1-H_2)\sqrt N\log_\mu\widehat\mu_N+o_p(1).
\end{equation}
The correction need not vanish even for balanced samples. Moreover, the earlier proof invokes an i.i.d. pooled-mean result although the pooled observations need not be identically distributed; Proposition~\ref{prop:pooled-clt} supplies the required linearization. The statistical cost of using the pooled covariance is quantified by the following general result, and in order to do so, we recall the definitions of the asymptotic errors.

\begin{Def}[\textbf{Actual asymptotic Type I and Type II errors}]\label{def:actual-errors}
The errors, calibrated at the nominal level $\alpha\in (0,1)$ are defined by:
\[
\boxed{
\alpha_{\mathrm{act}}
:= \lim_{N\to\infty}
{\mathbb{P}}_{H_0}
\left(T_{\mathrm{pool},N}^{2}>q_{d,1-\alpha}\right),\quad
\beta_{\mathrm{act}}(h)
:= \lim_{N\to\infty}
{\mathbb{P}}_{H_{1,N}(h)}
\left(T_{\mathrm{pool},N}^{2}\leq q_{d,1-\alpha}\right).
}
\]
Here $H_{1,N}(h)$ denotes a sequence of local alternatives
satisfying the limiting assumptions of Proposition~\ref{lem:pooled-size-distortion}.
\end{Def}

\begin{Prop}[\textbf{Cost of pooled-covariance misspecification and test calibration}]\label{lem:pooled-size-distortion}
{Suppose the baseline laws satisfy the assumptions of Proposition~\ref{prop:pooled-hotelling-limit} and $\Gamma^{(a)}\succ0$.} Put $\Sigma_{\lambda}:=\lambda \Sigma_1+(1-\lambda)\Sigma_2$, $\Gamma_{\mathrm{pool}}:=\Sigma_\lambda/[\lambda(1-\lambda)]$ and $\tau_N:=\tau_{\widehat\mu_N,\mu}$. Write $A:=\Gamma_{\mathrm{pool}}^{-1/2}\Gamma^{(a)}\Gamma_{\mathrm{pool}}^{-1/2}=Q\operatorname{diag}(\rho_1,\ldots,\rho_d)Q^\top$, where $Q$ is orthogonal, and let $q:=q_{d,1-\alpha}$ be the $(1-\alpha)$-quantile of $\chi_d^2$, with $0<\alpha<1$. For $h\in T_\mu M$, put $b:=Q^\top\Gamma_{\mathrm{pool}}^{-1/2}h$ and define
\begin{equation}\label{eq:misspecified-type-I}
\mathcal S_h(q):=\frac12+\frac1\pi\int_0^\infty\frac{\operatorname{Im}\left[e^{-itq}\prod_{j=1}^d\left[(1-2it\rho_j)^{-1/2}\exp\left(\frac{itb_j^2}{1-2it\rho_j}\right)\right]\right]}{t}\,dt.
\end{equation}
Then under $H_0$, the actual asymptotic Type~I error (Definition~\ref{def:actual-errors}) of the classical pooled Hotelling test, calibrated at nominal level $\alpha$ using $\chi_d^2$, is $\alpha_{\mathrm{act}}=\mathcal S_0(q)$, with explicit bounds
\[
1-F_{\chi_d^2}(q/\rho_{\min})\le\alpha_{\mathrm{act}}\le1-F_{\chi_d^2}(q/\rho_{\max}),
\]
where $\rho_{\min}:=\min_j\rho_j$ and $\rho_{\max}:=\max_j\rho_j$ are respectively the minimum and maximum eigenvalues of $A$ above. In particular, if $\Sigma_1=\Sigma_2\succ0$ and $H_1\neq H_2$, then $\alpha_{\mathrm{act}}>\alpha$. If $\Gamma^{(a)}=\rho\Gamma_{\mathrm{pool}}$, then $\alpha_{\mathrm{act}}=1-F_{\chi_d^2}(q/\rho)$. \revised{Under the local sampling regime of Theorem~\ref{thm:moving-clts}\textup{(ii)}, put $h:=H_2H_\lambda^{-1}a_1-H_1H_\lambda^{-1}a_2$. Whenever the perturbed laws violate $H_0$, the actual asymptotic Type~II error is $\beta_{\mathrm{act}}(h)=1-\mathcal S_h(q)$.}

Furthermore, let $\widehat\Gamma_{\mathrm{pool},N}:=N^2S_{p,N}/(nm)$ and suppose $\tau_N\widehat\Gamma_N^{(a)}\tau_N^*\inPrb\Gamma^{(a)}$ \revised{under $H_0$}. If $\widehat\rho_N$ consists of the eigenvalues of $\widehat\Gamma_{\mathrm{pool},N}^{-1/2}\widehat\Gamma_N^{(a)}\widehat\Gamma_{\mathrm{pool},N}^{-1/2}$, then the recalibrated $p$-value $p_N:=1-G(T_{\mathrm{pool},N}^2\mid\widehat\rho_N)$ satisfies $p_N\inD\mathrm{Unif}(0,1)$ under $H_0$, where $G(\cdot\mid\rho)$ is the weighted chi-square distribution function. The required inverses exist with probability tending to one.
\end{Prop}

\begin{proof}
\textbf{Under the null $H_0:$} By Proposition~\ref{prop:pooled-hotelling-limit}, the null limit under $H_0$ is the weighted chi-square distribution with weights $\rho_j$. For the specified local sampling regime, Theorem~\ref{thm:moving-clts}\textup{(ii)} gives $\sqrt N\,\tau_N\Delta_N^{(a)}\inD N_d(h,\Gamma^{(a)})$. 
The groupwise covariance consistency established in the proof of Lemma~\ref{lem:cov-consistency} gives $\tau_N\widehat\Sigma_{r,N}\tau_N^*\inPrb\Sigma_r$ under $H_0$, and hence $\tau_NS_{p,N}\tau_N^*\inPrb\Sigma_\lambda$. By contiguity, this convergence also holds under the specified local alternatives. Since $\widehat\Gamma_{\mathrm{pool},N}=N^2S_{p,N}/(nm)$ and $nm/N^2\to\lambda(1-\lambda)$, Slutsky's theorem gives $\tau_N\widehat\Gamma_{\mathrm{pool},N}\tau_N^*\inPrb\Sigma_\lambda/[\lambda(1-\lambda)]=\Gamma_{\mathrm{pool}}$.
{By invariance under $\tau_N$ and Slutsky's theorem, the pooled Hotelling statistic therefore converges to $Z^\top\Gamma_{\mathrm{pool}}^{-1}Z$, where $Z\sim N_d(h,\Gamma^{(a)})$.}

\textbf{Under the stated local alternatives:} \revised{With $Q$ as in the statement of this proposition, write $Y:=Q^\top\Gamma_{\mathrm{pool}}^{-1/2}Z$. Since $Z\sim N_d(h,\Gamma^{(a)})$ and $Q^\top A Q=\operatorname{diag}(\rho_1,\ldots,\rho_d)$, the coordinates of $Y$ are independent with $Y_j\sim N(b_j,\rho_j)$, where $b=(b_1,\ldots,b_d)^\top=Q^\top\Gamma_{\mathrm{pool}}^{-1/2}h$. Hence $Y_j\stackrel{\mathcal D}{=}\sqrt{\rho_j}G_j+b_j$, with independent $G_j\sim N(0,1)$. As $Q$ is orthogonal, $Z^\top\Gamma_{\mathrm{pool}}^{-1}Z=Y^\top Y$, so} the corresponding quadratic form is $\sum_{j=1}^d\rho_j(G_j+b_j/\sqrt{\rho_j})^2$, with independent $G_j\sim N(0,1)$.
 Its characteristic function is 
$$\prod_{j=1}^d\left[(1-2it\rho_j)^{-1/2}\exp\left(\frac{itb_j^2}{1-2it\rho_j}\right)\right],$$
obtained from \citet[equation~(2.3), p.~421]{Imhof1961} by taking one degree of freedom per component and noncentrality parameters $\delta_j^2=b_j^2/\rho_j$. Applying the characteristic-function inversion formula of \citet[Section~3, equations~(3.1)--(3.2), pp.~422--423]{Imhof1961} gives \eqref{eq:misspecified-type-I} and the asserted Type~I and Type~II errors. Specifically, the inversion yields the upper-tail probability $\mathcal S_h(q)$, so the Type~II error is its complement. \revised{Every $\rho_j>0$, because $\Gamma^{(a)}$ and $\Gamma_{\mathrm{pool}}$ are positive definite. The continuity-set form of the Portmanteau theorem for atom-free distributions therefore gives $\Prb_{0,N}(T_{\mathrm{pool},N}^2>q)\to\mathcal S_0(q)$ and $\Prb_{s,N}(T_{\mathrm{pool},N}^2\le q)\to1-\mathcal S_h(q)$, justifying the asserted rejection and non-rejection probabilities. Defining the statistic arbitrarily when an inverse fails does not affect these limits, since the probability of that event tends to zero.} The bounds follow from $\rho_{\min}\sum_jG_j^2\le\sum_j\rho_jG_j^2\le\rho_{\max}\sum_jG_j^2$. When $\Sigma_1=\Sigma_2\succ0$ and $H_1\neq H_2$, we have $\Gamma_{\mathrm{pool}}=\Gamma_{\mathrm{fix}}$, and Remark~\ref{rem:covariance-equality} implies $A\succeq I_d$ with $A\neq I_d$. Thus every $\rho_j\ge1$ and at least one $\rho_j>1$, giving strict inflation. Finally, \revised{the pooled-covariance consistency established above, together with the assumed consistency of $\widehat\Gamma_N^{(a)}$, implies by continuity of matrix inversion, positive-definite square roots and eigenvalues that} $\widehat\rho_N\inPrb(\rho_1,\ldots,\rho_d)$ up to ordering. \revised{Since $G(x\mid\rho)$ is jointly continuous in $x$ and the strictly positive weights $\rho$, Slutsky's theorem and the continuous mapping theorem give $p_N\inD1-G(W\mid\rho)$, where $W$ has distribution function $G(\cdot\mid\rho)$. The probability integral transform yields $1-G(W\mid\rho)\sim\mathrm{Unif}(0,1)$, proving the recalibrated $p$-value result.}
\end{proof}

\begin{Ex}[\textbf{Spherical size distortion}]\label{ex:huckemann-sphere}
On $S^2$, fix $\mu$ and an orthonormal basis $e_1,e_2$ of $T_\mu S^2$, and take $P_r=\frac14\sum_{j=1}^2(\delta_{\exp_\mu(t_re_j)}+\delta_{\exp_\mu(-t_re_j)})$, with $t_1=\pi/6$, $t_2=\pi/3$. Since $f(u)=\arccos(u)^2$ is strictly decreasing and convex on $[-1,1]$, Jensen's inequality gives $F_r(q)\ge f(\langle q,\mu\rangle\cos t_r)\ge t_r^2$, with equality only at $q=\mu$; thus $\mu$ is the unique global mean. The finite support lies in a common strongly convex ball of radius less than $\pi/2$. The spherical Hessian formula gives $H_r=h_rI_2$, $h_r=1+t_r\cot t_r>0$, and $\Sigma_r=s_rI_2$, $s_r=t_r^2/2$. For $n=m=N/2$, equation~\eqref{eq:arithmetic-two-group-covariance} yields $\Gamma_{\mathrm{fix}}=2(s_1+s_2)I_2\approx1.370778I_2$ and $\Gamma^{(a)}=8(h_2^2s_1+h_1^2s_2)(h_1+h_2)^{-2}I_2\approx1.522548I_2$. Hence
\begin{equation}\label{eq:huckemann-counterexample}
T_{\mathrm{pool},N}^2\inD\frac{1.522548\ldots}{1.370778\ldots}\chi_2^2\approx1.110717\,\chi_2^2,\qquad \alpha_{\mathrm{act}}(0.05)\approx0.0674.
\end{equation}
Thus balanced sampling does not ensure the claimed classical $\chi_2^2$ calibration for logarithmic coordinates, even though the regularity conditions hold. This does not concern other coordinate choices; the correctly studentized statistic of Corollary~\ref{cor:wald} retains its $\chi_d^2$ null limit.
\end{Ex}


\section{Concluding remarks and outlook}\label{sec:concluding-outlook}

Our results show that replacing the population basepoint by a pooled sample Fr\'echet mean has different effects on the two types of contrasts studied here. For contrasts of groupwise intrinsic sample means, the first-order basepoint terms cancel. For contrasts of arithmetic averages of individual logarithms, they need not cancel, and ignoring them can change the limiting distribution of a Hotelling-type test.

A similar question arises for orthogonal tangent projections on spheres. \citet[Section~3, Proposition~3.2]{MardiaPatrangenaru2005} use pooled sample mean directions as the poles for a two-sample Hotelling test and propose an $F$-approximation for concentrated data. Equality of population mean directions does not require equality of the lengths of the corresponding Euclidean mean vectors. When those lengths differ, estimating the pooled directions can contribute at first order to the projected contrast, an effect absent from a calculation that treats the poles as fixed.

The Procrustes tangent-space procedures of \citet[Chapter~7]{DrydenMardia1998} and \citet[Chapter~9, Section~9.1.2]{DrydenMardia2016} raise a related question. The latter explicitly describe choosing a tangent-space pole from the data and then treating it as fixed (p.~185). In their two-sample test, the pole is the pooled sample full Procrustes mean, and the observations are expressed in partial Procrustes tangent coordinates. Here, too, the effect of estimating the pole deserves examination. There is also a separate issue concerning the null hypothesis: equality of population full Procrustes mean shapes does not, in general, imply equality of the expected partial Procrustes tangent coordinates. A correction to the covariance alone would not resolve this issue.

A possible direction for further work is to derive moving-basepoint CLTs for these two types of tangent projections and investigate how to construct correctly calibrated Wald-type tests for the corresponding population hypotheses. The coordinate maps and pooled estimators differ from those considered here, so separate arguments would be required. In particular, the Procrustes case also requires careful consideration of the population hypothesis tested by the tangent-coordinate contrast. Our results for Riemannian logarithms do not establish these extensions. Furthermore, the cited Hotelling procedures are proposed as approximations for data with small variability; the present discussion does not claim that those approximations are generally invalid.

\appendix
\section{Geometric cancellation for intrinsic contrasts}\label{app:intrinsic-cancellation}
The estimate below establishes the basepoint replacement using only consistency, as used in Theorem~\ref{thm:moving-clts}\textup{(i)} and Remark~\ref{rem:smeary-moving}.
\begin{proof}[Proof of \eqref{eq:intrinsic-cancellation}]
Use $U=B_g(\mu,\rho)$ with $\overline{B_g(\mu,3\rho)}\subset B$
and the local smoothness properties specified in the proof of
Theorem~\ref{thm:moving-clts}\textup{(i)}. In the normal chart
$\phi:=\log_\mu$, its image is the Euclidean ball
$W=B(0,\rho)\subset T_\mu M$. All derivatives used below are bounded on
the relevant compact coordinate sets because the two-point logarithm and parallel transport are
smooth on a neighborhood of $\overline U\times\overline U$.
For $a,x,y\in U$, set
$G(a,x,y):=\tau_{a,\mu}(\log_a x-\log_a y)$.
For fixed $x$, the map $a\mapsto G(a,x,x)$ is identically zero:
its two logarithms are the same vector in $T_aM$ for every $a$.
Consequently, for every $\eta\in T_aM$ and every smooth curve
$a(t)$ in $U$ with $a(0)=a$, $\dot a(0)=\eta$,
\[
 D_aG(a,x,x)[\eta]
 =\left.\frac{d}{dt}\right|_{t=0}G(a(t),x,x)=0.
\]
The target space here is the fixed vector space $T_\mu M$, so this is an
ordinary derivative; no identification of varying tangent spaces is implicit.

To obtain a bound uniform in all three arguments, put
$\mathcal F(u,v):=\tau_{\exp_\mu u,\mu}
\log_{\exp_\mu u}(\exp_\mu v)$ and
$\widetilde G(u,v,w):=\mathcal F(u,v)-\mathcal F(u,w)$.
Bounded mixed derivatives on $\overline W\times\overline W$ and the
fundamental theorem of calculus give
$$
 D_1\widetilde G(u,v,w) = \int_0^1 D_2D_1\mathcal F(u,w+t(v-w))[v-w]\,dt, \qquad
 \|D_1\widetilde G(u,v,w)\|_{\mathrm{op}} \le K\|v-w\|.
$$

The segment $w+t(v-w)$ lies in the Euclidean ball $W$.
To compare coordinate and Riemannian distances, take any minimizing geodesic $\gamma$ from $x$ to $y$.
Each $z\in\gamma$ satisfies
$d_g(\mu,z)\le d_g(\mu,x)+d_g(x,z)\le d_g(\mu,x)+d_g(x,y)<3\rho$,
so $\gamma$ lies in $B_g(\mu,3\rho)\subset B$.
The derivative of $\phi$ is bounded on $\overline{B_g(\mu,3\rho)}$;
integration along $\gamma$ therefore gives
$\|\phi(x)-\phi(y)\|\le K_1d_g(x,y)$.
Integrating the preceding estimate along the segment from $0$ to
$u=\phi(a)$, and using $\|u\|=d_g(a,\mu)$, yields
\begin{equation}\label{eq:intrinsic-cancellation}
 \|G(a,x,y)-G(\mu,x,y)\|
 \le C d_g(a,\mu)d_g(x,y),\qquad a,x,y\in U.
\end{equation}

\end{proof}

\section{Proof of the empirical moving-logarithm lemma}\label{app:empirical-moving-log}
\begin{proof}[Proof of Lemma~\ref{lem:empirical-moving-log}]
\revised{We first work under baseline sampling; $G_r$ and $H_r$ below are determined by $P_r$.}

\emph{Case (A).} Put $\xi_N:=\log_\mu a_N$. Compact support separated from $\Cut(\mu)$ supplies the uniform smooth product neighborhood used in \eqref{eq:common-hessian-expansion}. Thus $h_a(x)=h_\mu(x)+L_x\log_\mu a+O(d_g(a,\mu)^2)$ uniformly for $x\in C$, with $L_x=-K_\mu(x)/2$, 
where $K_\mu(x):=\Hess_{p=\mu}^{\sharp}d_g(p,x)^2$ denotes the $(1,1)$-Hessian operator of the squared Riemannian distance with respect to its first argument, as defined immediately before Equation \eqref{eq:common-hessian-expansion}.
Compactness and smoothness justify differentiation under the integral, giving $P_rL=-H_r/2$. Consequently $P_{r,n_r}h_{a_N}=U_{r,N}+(P_{r,n_r}L)\xi_N+O_p(N^{-1})$. The law of large numbers gives $P_{r,n_r}L\inPrb-H_r/2$; multiplication by $\sqrt N\xi_N=O_p(1)$ proves \eqref{eq:empirical-moving-log}.

\emph{Case (B).} By \citet[Theorem~1, equations~(13)--(15)]{HotzLeWood2024} and $G_r(\mu)=0$, the population logarithm field satisfies
\begin{equation}\label{eq:population-moving-log}
\tau_{a,\mu}G_r(a)=-\frac12H_r\log_\mu a+o(d_g(a,\mu)),\qquad a\to\mu.
\end{equation}
Their equations~(14)--(15) give $-H_r/2=\int_M L_x\,dP_r(x)+J_{P_r}(\mu)$, where $L_x$ is the pointwise covariant derivative of $p\mapsto\log_p x$ at $\mu$ for $x\notin\Cut(\mu)$ and is set to zero on the $P_r$-null set $\Cut(\mu)$. The cut-locus term $J_{P_r}(\mu)$ is defined in \citet[equation~(20)]{HotzLeWood2024} and obtained using \citet[Proposition~3]{HotzLeWood2024}; it is therefore included in $H_r$.

Since $P_rh_a=\tau_{a,\mu}G_r(a)$ and $P_rh_\mu=G_r(\mu)=0$, adding and subtracting population integrals gives the exact identity
\begin{equation}\label{eq:empirical-moving-decomposition}
P_{r,n_r}h_{a_N}=P_{r,n_r}h_\mu+\tau_{a_N,\mu}G_r(a_N)+(P_{r,n_r}-P_r)(h_{a_N}-h_\mu).
\end{equation}
For fixed $\delta,\varepsilon>0$, the probability that the norm of the final term, multiplied by $\sqrt{n_r}$, exceeds $\varepsilon$ is bounded by $\Prb(a_N\notin B_\delta(\mu))$ plus the probability inside~\eqref{eq:moving-log-equicontinuity}. Taking $\limsup_{N\to\infty}$ and then letting $\delta\downarrow0$, consistency of $a_N$ and~\eqref{eq:moving-log-equicontinuity} show that this term is $o_p(n_r^{-1/2})=o_p(N^{-1/2})$. The remainder in~\eqref{eq:population-moving-log}, evaluated at $a_N$, is likewise $o_p(N^{-1/2})$ by consistency and $\sqrt N\log_\mu a_N=O_p(1)$. Substitution into~\eqref{eq:empirical-moving-decomposition} proves~\eqref{eq:empirical-moving-log}. The result holds simultaneously because $k$ is fixed.

\revised{For local sampling, regard the same measurable sequence $a_N$ as a function on the baseline sample space. Mutual contiguity, proved in Section~\ref{sec:moving-clts}, transfers its consistency and root-$N$ tightness from $\Prb_{s,N}$ to $\Prb_{0,N}$. Applying the preceding baseline proof in case \textup{(A)} or \textup{(B)} makes the remainder in \eqref{eq:empirical-moving-log} $o_{\Prb_{0,N}}(N^{-1/2})$; contiguity in the forward direction makes the same remainder $o_{\Prb_{s,N}}(N^{-1/2})$. Throughout case \textup{(B)}, $G_r(\mu)=P_rh_\mu=0$ is the baseline identity, whereas $P_{r,N}h_\mu=a_r/\sqrt N$ generally does not vanish. The exact decomposition \eqref{eq:empirical-moving-decomposition} still uses $P_r$ and $G_r$, not perturbed population fields; its stochastic remainder is transferred as a whole. Thus the Hessian remains the full baseline $H_r$, including any cut-locus contribution, without requiring the perturbed laws themselves to satisfy the smoothness assumptions of \citet{HotzLeWood2024}.}
\end{proof}

\begin{Remark}[\textbf{Empirical-process condition}]\label{rem:cut-locus-equicontinuity}
Theorem~1 of \citet{HotzLeWood2024} yields the population expansion~\eqref{eq:population-moving-log}, but does not itself establish~\eqref{eq:moving-log-equicontinuity}. Step~2 of their proof of Theorem~2 (pp.~1240--1242, equations~(43)--(48)) invokes a Donsker argument for the corresponding i.i.d.\ empirical process. In the printed argument, however, the $O(d_g(x,y))$ volume bound for the fixed set $A_\delta(\mu)$ used in equation~(46) is not justified by assumption~\textup{(iv)}, which gives $\operatorname{vol}_M(A_t(\mu))=O(t)$ as $t\downarrow0$; moreover, the asserted bracketing-entropy conclusion on p.~1242 is not established by the $L^2$ increment bound~(47) alone. Therefore~\eqref{eq:moving-log-equicontinuity} is retained as an explicit assumption rather than inferred from \textup{(iii)--(v)}. The $M$-estimator theorem of \citet[Theorem~5.23]{vanDerVaart1998} does not itself supply this empirical moving-logarithm equicontinuity.
\end{Remark}

\section{Proof of the pooled Hotelling limit}\label{app:pooled-hotelling}

\begin{Prop}[\textbf{Limiting law of the pooled Hotelling statistic}]\label{prop:pooled-hotelling-limit}
Under $H_0$, Assumption~\ref{as:compact-support}, and the hypotheses of Theorem~\ref{thm:moving-clts}\textup{(ii)}, suppose $\Sigma_\lambda:=\lambda\Sigma_1+(1-\lambda)\Sigma_2$ is positive definite. Then, for arbitrary $\lambda\in(0,1)$,
\begin{equation}\label{eq:pooled-hotelling-general-limit}
T_{\mathrm{pool},N}^2:=\frac{nm}{N}\left\langle
 V_{1,N}^{(a)}-V_{2,N}^{(a)},
 S_{p,N}^{-1}(V_{1,N}^{(a)}-V_{2,N}^{(a)})
 \right\rangle\inD\sum_{j=1}^d\rho_j\chi_{1,j}^2,
\end{equation}
where $\chi_{1,j}^2$ are independent $\chi_1^2$ variables and $\rho_j$ are the eigenvalues of $A:=\lambda(1-\lambda)\Sigma_\lambda^{-1/2}\Gamma^{(a)}\Sigma_\lambda^{-1/2}$. The limit is $\chi_d^2$ if and only if $\Gamma^{(a)}=\Sigma_\lambda/[\lambda(1-\lambda)]$.
\end{Prop}
\begin{proof}
Put $\tau_N:=\tau_{\widehat\mu_N,\mu}$ and
$Z_N:=\sqrt N\,\tau_N\Delta_N^{(a)}$.
Theorem~\ref{thm:moving-clts}\textup{(ii)} gives
$Z_N\inD Z\sim N_d(0,\Gamma^{(a)})$.
By the definition of the within-group pooled sample covariance,
\[
S_{p,N}
=\frac{(n-1)\widehat\Sigma_{1,N}
+(m-1)\widehat\Sigma_{2,N}}{N-2}.
\]
The groupwise covariance consistency established in the proof of
Lemma~\ref{lem:cov-consistency} therefore yields
\[
\tau_NS_{p,N}\tau_N^*
\inPrb\lambda\Sigma_1+(1-\lambda)\Sigma_2
=:\Sigma_\lambda.
\]

Since $\tau_N$ is a linear isometry, we may rewrite the statistic as
\[
T_{\mathrm{pool},N}^2
=\frac{nm}{N^2}\left\langle Z_N,
(\tau_NS_{p,N}\tau_N^*)^{-1}Z_N\right\rangle.
\]
Now $\tau_NS_{p,N}\tau_N^*\inPrb\Sigma_\lambda$ and
$nm/N^2\to\lambda(1-\lambda)$. Since $\Sigma_\lambda$ is positive definite, the required inverses exist with probability tending to one. Hence, by Slutsky's theorem and
continuity of matrix inversion,
\[
T_{\mathrm{pool},N}^2\inD
\lambda(1-\lambda)Z^\top\Sigma_\lambda^{-1}Z.
\]
Put $W:=\sqrt{\lambda(1-\lambda)}\,\Sigma_\lambda^{-1/2}Z$, so that
$W\sim N_d(0,A)$. By the spectral theorem,
\[
A=Q\operatorname{diag}(\rho_1,\ldots,\rho_d)Q^\top
\]
for some orthogonal matrix $Q$. Consequently, the coordinates
of $Q^\top W$ are independent centered Gaussian variables
with variances $\rho_1,\ldots,\rho_d$, and therefore
\[
\lambda(1-\lambda)Z^\top\Sigma_\lambda^{-1}Z
=W^\top W
\overset{\mathcal D}{=}
\sum_{j=1}^d\rho_jG_j^2,
\]
where $G_1,\ldots,G_d$ are independent $N(0,1)$ variables.

Finally, if this distribution equals $\chi_d^2$, comparison
of their means and variances gives
$\sum_{j=1}^d\rho_j=d$ and
$\sum_{j=1}^d\rho_j^2=d$, respectively, since
$\mathbb E[G_j^2]=1$ and $\operatorname{Var}(G_j^2)=2$.
Thus
$\sum_{j=1}^d(\rho_j-1)^2=0$, implying
$\rho_j=1$ for every $j$. Equivalently, $A=I_d$, or
$\Gamma^{(a)}=\Sigma_\lambda/[\lambda(1-\lambda)]$. The converse is immediate.
\end{proof}

\begin{Cor}[\textbf{Balanced samples or equal group covariances}]\label{cor:pooled-hotelling-huckemann}
Under the assumptions of Proposition~\ref{prop:pooled-hotelling-limit}, suppose additionally that either $\lambda=1/2$ or $\Sigma_1=\Sigma_2$, as in the two-sample setting of \citet[Section~6.1, Theorem~10]{Huckemann2012}. Then
\[
T_{\mathrm{pool},N}^2\inD\sum_{j=1}^d\rho_j\chi_{1,j}^2,
\]
where $\rho_j$ are now the eigenvalues of $A:=\Gamma_{\mathrm{fix}}^{-1/2}\Gamma^{(a)}\Gamma_{\mathrm{fix}}^{-1/2}$, with $\Gamma_{\mathrm{fix}}:=\Sigma_1/\lambda+\Sigma_2/(1-\lambda)$. In particular, this limiting distribution equals $\chi_d^2$ if and only if $\Gamma^{(a)}=\Gamma_{\mathrm{fix}}$.
\end{Cor}
\begin{proof}
Since
\(
\Sigma_\lambda-\lambda(1-\lambda)\Gamma_{\mathrm{fix}}
=(2\lambda-1)(\Sigma_1-\Sigma_2)=0
\)
under either of the stated conditions, we have $\Sigma_\lambda=\lambda(1-\lambda)\Gamma_{\mathrm{fix}}$ under these conditions. Consequently,
\(
\lambda(1-\lambda)\Sigma_\lambda^{-1/2}
\Gamma^{(a)}\Sigma_\lambda^{-1/2}
=\Gamma_{\mathrm{fix}}^{-1/2}
\Gamma^{(a)}\Gamma_{\mathrm{fix}}^{-1/2}.
\)
The result follows immediately from Proposition~\ref{prop:pooled-hotelling-limit}.
\end{proof}

\section{Table of symbols and their meanings}

\begin{center}
\small
\setlength{\LTleft}{0pt}
\setlength{\LTright}{0pt}
\renewcommand{\arraystretch}{1.15}
\begin{longtable}{p{0.24\textwidth}p{0.70\textwidth}}
\caption{Table of notation.}\label{tab:notation}\\
\hline
\textbf{Symbol} & \textbf{Meaning} \\
\hline
\endfirsthead
\hline
\textbf{Symbol} & \textbf{Meaning} \\
\hline
\endhead
\hline
\endfoot

$M$ & Complete connected $d$-dimensional Riemannian manifold. \\

$d_g(\cdot,\cdot)$ & Riemannian distance on $M$. \\

$\mu$ & Common population intrinsic Fréchet mean under $H_0$. \\

$\operatorname{Cut}(a)$ & Cut locus of $a\in M$. \\

$T_aM$ & Tangent space of $M$ at $a\in M$. \\

$\log_a x$ & Riemannian logarithm of $x$ at basepoint $a$, with the zero-extension convention on the cut locus when explicitly adopted in the paper. \\

$\exp_\mu$ & Riemannian exponential map at $\mu$. \\

$\tau_{a,\mu}$ & Parallel transport from $T_aM$ to $T_\mu M$ along the unique short minimizing geodesic from $a$ to $\mu$. \\

$\tau_N$ & Shorthand for $\tau_{\widehat\mu_N,\mu}$. \\

$X_j^{(r)}$ & $j$th observation in group $r$. \\

$P_r$ & Probability law of group $r$. \\

$n_r$ & Sample size of group $r$. \\

$N$ & Total sample size, $N=\sum_{r=1}^k n_r$. \\

$\lambda_{r,N}$ & Empirical sample proportion $n_r/N$. \\

$\lambda_r$ & Asymptotic sample proportion, $\lambda_{r,N}\to\lambda_r\in(0,1)$. \\

$\lambda$ & In the two-group case, $\lambda:=\lambda_1$ and $1-\lambda=\lambda_2$. \\

$F_r(p)$ & Population Fréchet function for group $r$, $F_r(p)=E[d_g(p,X^{(r)})^2]$. \\

$F_{r,n_r}(p)$ & Empirical Fréchet function for group $r$, $F_{r,n_r}(p)=n_r^{-1}\sum_{j=1}^{n_r} d_g(p,X_j^{(r)})^2$. \\

$F_N$ & Pooled empirical Fréchet function, $F_N=\sum_{r=1}^k \lambda_{r,N}F_{r,n_r}$. \\

$F_\lambda$ & Affine population combination $F_\lambda=\sum_{r=1}^k \lambda_r F_r$. \\

$\widehat\mu_{r,n_r}$ & Groupwise sample intrinsic Fréchet mean for group $r$. \\

$\widehat\mu_N$ & Pooled sample intrinsic Fréchet mean. \\

$\xi_N$ & Pooled mean log-coordinate, $\xi_N:=\log_\mu\widehat\mu_N$. \\

$H_r$ & Population Hessian operator of $F_r$ at $\mu$, i.e. $H_r=\operatorname{Hess}_\mu^\sharp F_r$. \\

$H_\lambda$ & Affine Hessian combination $H_\lambda=\sum_{r=1}^k \lambda_r H_r$; in the two-group case, $H_\lambda=\lambda H_1+(1-\lambda)H_2$. \\

$\Sigma_r$ & Groupwise tangential covariance operator, $\Sigma_r=\operatorname{Cov}(\log_\mu X^{(r)})$. \\

$\Sigma_\lambda$ & Affine covariance combination $\Sigma_\lambda=\sum_{r=1}^k \lambda_r\Sigma_r$; in the two-group case, $\Sigma_\lambda=\lambda\Sigma_1+(1-\lambda)\Sigma_2$. \\

$U_{r,N}$ & Groupwise arithmetic average of fixed-basepoint logarithms, $U_{r,N}=n_r^{-1}\sum_{j=1}^{n_r}\log_\mu X_j^{(r)}$. \\

$V_{r,N}^{(a)}$ & Groupwise arithmetic average of moving-basepoint logarithms, $V_{r,N}^{(a)}=n_r^{-1}\sum_{j=1}^{n_r}\log_{\widehat\mu_N}X_j^{(r)}$. \\

$\Delta_{\mu,N}^{(i)}$ & Fixed-basepoint intrinsic contrast vector formed from the groupwise sample Fréchet means. \\

$\Delta_N^{(i)}$ & Moving-basepoint intrinsic contrast vector formed from the groupwise sample Fréchet means. \\

$\Delta_N^{(a)}$ & Moving-basepoint arithmetic contrast vector; in the two-group case, $\Delta_N^{(a)}=V_{1,N}^{(a)}-V_{2,N}^{(a)}$. \\

$Q_r$ & Groupwise intrinsic CLT covariance, $Q_r=4H_r^{-1}\Sigma_r H_r^{-1}$. \\

$\Gamma^{(i)}$ & Limiting covariance operator of the intrinsic moving-basepoint contrast CLT. \\

$\Gamma^{(a)}$ & Limiting covariance operator of the arithmetic moving-basepoint contrast CLT. \\

$B_{rs}$ & Coefficients in the arithmetic contrast linearization, $B_{rs}=(\delta_{rs}-\delta_{ks})I_d-\lambda_s(H_r-H_k)H_\lambda^{-1}$. \\

$B$ & Block matrix $(B_{rs})_{r<k,\;s\le k}$ used in the representation $\Gamma^{(a)}=BD_\lambda B^*$. \\

$D_\lambda$ & Block-diagonal operator $\operatorname{diag}(\Sigma_1/\lambda_1,\ldots,\Sigma_k/\lambda_k)$. \\

$Z_{rj,N}$ & Moving-basepoint log-observation, $Z_{rj,N}:=\log_{\widehat\mu_N}X_j^{(r)}$. \\

$\overline Z_{r,N}$ & Mean of the $Z_{rj,N}$, i.e. $\overline Z_{r,N}=n_r^{-1}\sum_{j=1}^{n_r} Z_{rj,N}$. \\

$\widehat\Sigma_{r,N}$ & Sample covariance operator (or matrix in a chosen orthonormal basis) of the moving-basepoint log-observations in group $r$. \\

$\widehat H_{r,N}$ & Empirical Hessian operator at $\widehat\mu_N$, defined from the sample Fréchet function of group $r$. \\

$\widehat H_{\lambda,N}$ & Affine empirical Hessian combination, $\widehat H_{\lambda,N}=\sum_{r=1}^k \lambda_{r,N}\widehat H_{r,N}$. \\

$\widehat\Gamma_N^{(i)}$ & Plug-in estimator of $\Gamma^{(i)}$. \\

$\widehat\Gamma_N^{(a)}$ & Plug-in estimator of $\Gamma^{(a)}$. \\

$T_N^{(i)}$ & Corrected Wald statistic for the intrinsic moving-basepoint CLT. \\

$T_N^{(a)}$ & Corrected Wald statistic for the arithmetic moving-basepoint CLT. \\

$S_{p,N}$ & Ordinary within-group pooled sample covariance matrix of the coordinates $\log_{\widehat\mu_N}X_j^{(r)}$ in the two-group case. \\

$\Gamma_{\mathrm{pool}}$ & Population covariance targeted by the rescaled pooled covariance, $\Gamma_{\mathrm{pool}}:=\Sigma_\lambda/[\lambda(1-\lambda)]$. \\

$\widehat\Gamma_{\mathrm{pool},N}$ & Rescaled pooled covariance estimator, $\widehat\Gamma_{\mathrm{pool},N}:=N^2S_{p,N}/(nm)$. \\

$\Gamma_{\mathrm{fix}}$ & Fixed-basepoint arithmetic covariance in the two-group case, $\Gamma_{\mathrm{fix}}=\Sigma_1/\lambda+\Sigma_2/(1-\lambda)$. \\

$T_{\mathrm{pool},N}^2$ & Classical pooled Hotelling-type statistic based on the arithmetic contrast and pooled covariance normalization. \\

$D$ & In Remark~\ref{rem:covariance-equality}, $D:=H_1-H_2$. \\

$C$ & In Remark~\ref{rem:covariance-equality}, $C:=DH_\lambda^{-1}=(H_1-H_2)H_\lambda^{-1}$. \\

$A$ & Misspecification matrix comparing the correct and pooled covariances, $A=\Gamma_{\mathrm{pool}}^{-1/2}\Gamma^{(a)}\Gamma_{\mathrm{pool}}^{-1/2}$; under $\lambda=1/2$ or $\Sigma_1=\Sigma_2$, this equals $\Gamma_{\mathrm{fix}}^{-1/2}\Gamma^{(a)}\Gamma_{\mathrm{fix}}^{-1/2}$. \\

$\rho_1,\ldots,\rho_d$ & Eigenvalues of the matrix/operator $A$. \\

$q_{d,1-\alpha}$ & $(1-\alpha)$-quantile of the $\chi_d^2$ distribution. \\

$\mathcal S_h(q)$ & Imhof-type upper-tail probability formula used for the asymptotic Type I and Type II error calculations under covariance misspecification. \\

$\alpha_{\mathrm{act}}$ & Actual asymptotic Type I error of the classical pooled Hotelling test under misspecified covariance normalization. \\

$\beta_{\mathrm{act}}$ & Actual asymptotic Type II error under local alternatives. \\

$\widehat\rho_N$ & Estimated eigenvalues used in the recalibrated weighted-chi-square $p$-value. \\

$G(\,\cdot\mid\rho)$ & Distribution function of the weighted chi-square law determined by the weights $\rho=(\rho_1,\ldots,\rho_d)$. \\

$u,v,h$ & Generic tangent vectors, usually belonging to $T_\mu M$ unless another tangent space is specified. \\

$I_d$ & Identity operator (or identity matrix, once an orthonormal basis is fixed) on a $d$-dimensional tangent space. \\

$\delta_{rs}$ & Kronecker delta. \\

\end{longtable}
\end{center}

\section*{Acknowledgments and Funding}

I am grateful to the Hausdorff Research Institute for Mathematics
(HIM) for its hospitality during my invited stay there in the Fall of 2026, and to Stephan Huckemann for suggesting the inclusion of the smeary case.

\section*{Declaration of Generative AI Use}
The author used ChatGPT (GPT-6, OpenAI) to facilitate certain 
discussions of mathematical ideas and to assist with
manuscript drafting and literature searches. The author takes full responsibility for the final content.

\bibliographystyle{plainnat}
\bibliography{References}

\end{document}